\documentclass[11pt]{article}
\pdfoutput=1
\usepackage{amssymb,latexsym,amsmath}
\usepackage{epsfig}
\usepackage{caption}
\usepackage{graphicx}
\def\mindex#1{\index{#1}}

\def\sq{\hbox{\rlap{$\sqcap$}$\sqcup$}}
\def\qed{\ifmmode\sq\else{\unskip\nobreak\hfil
\penalty50\hskip1em\null\nobreak\hfil\sq
\parfillskip=0pt\finalhyphendemerits=0\endgraf}\fi\medskip}

\long\def\defbox#1{\framebox[.9\hsize][c]{\parbox{.85\hsize}{%
\parindent=0pt
\baselineskip=12pt plus .1pt      % STYLE
\parskip=6pt plus 1.5pt minus 1pt % CHANGES
 #1}}}

\long\def\beginbox#1\endbox{\subsection*{}%
\hbox{\hspace{.05\hsize}\defbox{\medskip#1\bigskip}}%
\subsection*{}}

\def\endbox{}

\newsavebox{\junk}
\savebox{\junk}[1.6mm]{\hbox{$|\!|\!|$}}

\def\limsup{\mathop{\rm lim\ sup}}
\def\liminf{\mathop{\rm lim\ inf}}

\newcommand{\field}[1]{\mathbb{#1}}

\def\Re{\field{R}}

\def\bfmath#1{{\mathchoice{\mbox{\boldmath$#1$}}%
{\mbox{\boldmath$#1$}}%
{\mbox{\boldmath$\scriptstyle#1$}}%
{\mbox{\boldmath$\scriptscriptstyle#1$}}}}

\def\bfmX{\bfmath{X}}

\def\bfmY{\bfmath{Y}}

\def\bfmhhaY{\bfmath{\hhaY}} %\widehat{\widehat{Y}}}}
\def\bfmhhaY{\hbox to 0pt{$\widehat{\bfmY}$\hss}\widehat{\phantom{\raise 1.25pt\hbox{$\bfmY$}}}}

\def\til={{\widetilde =}}

\def\clF{{\cal F}}

\def\clM{{\cal M}}

 \def\FRAC#1#2#3{\genfrac{}{}{}{#1}{#2}{#3}}

\def\ddtp{{\mathchoice{\FRAC{1}{d^{\hbox to 2pt{\rm\tiny +\hss}}}{dt}}%
{\FRAC{1}{d^{\hbox to 2pt{\rm\tiny +\hss}}}{dt}}%
{\FRAC{3}{d^{\hbox to 2pt{\rm\tiny +\hss}}}{dt}}%
{\FRAC{3}{d^{\hbox to 2pt{\rm\tiny +\hss}}}{dt}}}}

\def\half{{\mathchoice{\FRAC{1}{1}{2}}%
{\FRAC{1}{1}{2}}%
{\FRAC{3}{1}{2}}%
{\FRAC{3}{1}{2}}}}

\def\eqdef{\mathbin{:=}}

\def\Prob{P}
\def\Expect{E}

\def\average#1,#2,{{1\over #2} \sum_{#1}^{#2}}

\def\eye(#1){{\bf(#1)}\quad}

\newtheorem{theorem}{Theorem}[section]

\def\eq#1/{(\ref{e:#1})}

\newcommand{\beqn}[1]{\notes{#1}%
\begin{eqnarray} \elabel{#1}}

\newcommand{\eeqn}{\end{eqnarray} }

\newcommand{\beq}[1]{\notes{#1}%
\begin{equation}\elabel{#1}}

\newcommand{\eeq}{\end{equation}}

\def\bdes{\begin{description}}
\def\edes{\end{description}}

\newcounter{rmnum}
\newenvironment{romannum}{\begin{list}{{\upshape (\roman{rmnum})}}{\usecounter{rmnum}
\setlength{\leftmargin}{14pt}
\setlength{\rightmargin}{8pt}
\setlength{\itemindent}{-1pt}
}}{\end{list}}

\newcounter{anum}

{\end{list}}

\def\ass(#1:#2){(#1\ref{#1:#2})}

\def\ritem#1{
\item[{\sf \ass(\current_model:#1)}]
}

\newenvironment{recall-ass}[1]{%
\begin{description}
\def\current_model{#1}}{
\end{description}
}

\newcommand{\bd}{\begin{description}}
\newcommand{\ed}{\end{description}}
\newcommand{\bt}{\begin{theorem}}
\newcommand{\et}{\end{theorem}}
\newcommand{\ba}{\begin{array}{rcl}}
\newcommand{\ea}{\end{array}}

\def\stp{{\cal T}}

\newlength{\noteWidth}
\long\def\notes#1{\ifinner
           {\tiny #1}
           \else
           \marginpar{\parbox[t]{\noteWidth}{\raggedright\tiny #1}}
       \fi\typeout{#1}}

       \newtheorem{thm}{\bf{Theorem}}[section]
       \newtheorem{cor}{\bf{Corollary}}[section]
       \newtheorem{lem}{\bf{Lemma}}[section]

       \newtheorem{defn}{\bf{Definition}}[section]
  \newtheorem{remark}{Remark}
       \newtheorem{assumption}{\bf{Assumption}}[section]

\allowdisplaybreaks

\def\qed{\hfill $\diamond$}
\fussy
\begin{document}
\title{Stationary and Ergodic Properties of Stochastic Non-Linear Systems Controlled over Communication Channels}
\author{Serdar Y\"uksel
\thanks{Department of Mathematics and
    Statistics, Queen's University, Kingston, Ontario, Canada, K7L
    3N6.  Email: yuksel@mast.queensu.ca. This research was
    partially supported by the Natural Sciences and Engineering
    Research Council of Canada (NSERC). Some of the results in this paper have been presented at the 2016 IEEE International Symposium on Information Theory (ISIT).}
}
%\begin{document}

\maketitle
\begin{abstract}
This paper is concerned with the following problem: Given a stochastic non-linear system controlled over a noisy channel, what is the largest class of channels for which there exist coding and control policies so that the closed loop system is stochastically stable? Stochastic stability notions considered are stationarity, ergodicity or asymptotic mean stationarity. We do not restrict the state space to be compact, for example systems considered can be driven by unbounded noise. Necessary and sufficient conditions are obtained for a large class of systems and channels. A generalization of Bode's Integral Formula for a large class of non-linear systems and information channels is obtained. The findings generalize existing results for linear systems. 
\end{abstract}

%\begin{keyword}
%{\bf Key words}: Networked control, stochastic stability, ergodicity, non-linear systems
%%\kwd{\LaTeXe}
%\end{keyword}
%
%\begin{keyword} {\bf AMS subject classifications}: 93E03, 60J20
%\end{keyword}
%
\section{Introduction}

Consider an $N$-dimensional controlled non-linear system described by the discrete-time equations
\begin{eqnarray}
x_{t+1} &=& f(x_t,u_t, w_t), \label{systemMODEL1} 
%y_t&=& g(x_t,v_t), t \geq 0, \label{systemMODEL}
\end{eqnarray}
for a (Borel measurable) function $f$, with $\{w_t\}$ being an independent and identically distributed (i.i.d) system noise process with $w_t \sim \nu$.

\begin{figure}[h]
        \begin{center}
        \includegraphics[height=2.5cm,width=8cm]{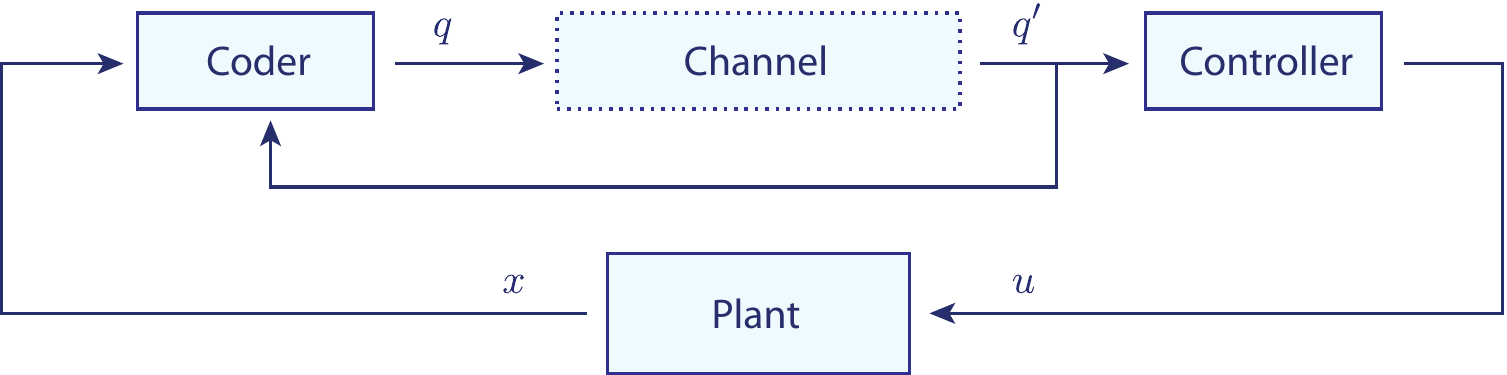}
\caption{Control of a system over a noisy channel.\label{LLL1}}
        \end{center}
\end{figure}

This system is connected over a noisy channel with a finite capacity to a controller, as shown in Figure~\ref{LLL1}. The controller has access to the information it has received through the channel. A source coder maps the source symbols, state values, to corresponding channel inputs. The channel inputs are transmitted through a channel; we assume that the channel is a finite alphabet channel with input alphabet ${\cal M}$ and output alphabet ${\cal M}'$.

We refer by a {\em Coding Policy} $\Pi$, a sequence of functions $\{\gamma^e_t, t \geq 0 \}$ which are causal such that the channel input at time $t$, $q_t \in {\cal M}$, under $\Pi^{comp}$ is generated by a function of its local information, that is, \[q_t=\gamma^e_t({\cal I}^e_t),\] where ${\cal I}^e_t=\{x_{[0,t]}, q'_{[0,t-1]}\}$ and $q_t \in {\cal M}$, the channel input alphabet given by ${\cal M} := \{1,2,\dots, M \},$ for $0 \leq t \leq T-1$. Here, we have the notation for $t \geq 1$: $x_{[0,t-1]} = \{x_s, 0 \leq s \leq t-1 \}.$

The channel maps $q_t$ to $q'_t$ in a stochastic fashion so that $P(q'_t|q_t, q_{[0,t-1]}, q'_{[0,t-1]})$ is a conditional probability measure on ${\cal M}'$ for all $t \in \mathbb{Z}_+$. If this expression is equal to $P(q'_t|q_t)$, the channel is said to be a memoryless channel, that is, the past variables do not affect the channel output $q'_t$ given the current channel input $q_t$. Even though in this paper we will consider discrete alphabet channels, the analysis is also applicable to a large class of continuous alphabet channels (through an appropriate quantized approximation of the channel; see e.g. \cite{ElGamalKim}). 

The receiver/controller, upon receiving the information from the channel, generates its decision at time $t$, also causally: An admissible causal controller policy is a sequence of functions $\gamma=\{\gamma_t\}$ such that $$\gamma_t : {\cal M}'^{t+1} \to \mathbb{R}^m, \quad \quad t \geq 0,$$
so that $u_t = \gamma_t(q'_{[0,t]})$. We call such encoding and control policies, {\em causal} or {\em admissible}.

In the networked control literature, the goal in the encoder/controller design is typically either to optimize the system according to some performance criterion or stabilize the system. For stabilization, linear systems have been studied extensively where the goal has been to identify conditions so that the controlled state is stochastically stable, as we review briefly later.

This paper is concerned with necessary and sufficient conditions on information channels in a networked control system for which there exist coding and control policies such that the controlled system is stochastically stable in one or more of the following senses: (i) The state $\{x_t\}$ and the coding and control parameters lead to a stable (positive Harris recurrent) Markov chain and (ii) $\{x_t\}$ is asymptotically stationary, or asymptotically mean stationary (AMS) and satisfies Birkhoff's sample path ergodic theorem (see Section \ref{ergodicTheory} for a review of these concepts), (iii) $\{x_t\}$ is ergodic.

In the remainder of this section, we will be providing a literature review, first for non-linear systems and then briefly for linear systems in the context of the goals of this paper and highlight the contributions of the paper. Section \ref{Bode} develops some supporting results and a generalization of Bode's Integral Formula for non-linear systems and general information channels. Section \ref{Ergodicity} develops conditions for ergodicity and asymptotic mean stationarity of the controlled system. Section \ref{Stationary} establishes conditions for stationarity of the controlled system under structured (stationary) coding and control policies. Section \ref{Sec3} presents an ergodic construction for a non-linear system driven by additive Gaussian noise and controlled over discrete noiseless channels. Section \ref{Sec4} contains some concluding remarks.

\subsection{Some notation and preliminaries}
Let $x$ be an $\mathbb{X}-$valued random variable, where $\mathbb{X}$ is countable. The {\em entropy} of $x$ is defined as $H(x) = -\sum_{z \in \mathbb{X}} p(z) \log_2(p(z))\,,$ where $p$ is the probability mass function (pmf) of the random variable $x$. If $x$ is an $\mathbb{R}^n-$valued random variable, and the probability measure induced by $x$ is absolutely continuous with respect to the Lebesgue measure, the {\em (differential) entropy} of $x$ is defined by $h(x) = -\int_{\mathbb{X}} p(x) \log_2(p(x))dx\,,$ where $p(\cdot)$ is the probability density function (pdf) of $x$. 

The {\em Mutual Information} between a discrete (continuous) random variable $x$, and another discrete (continuous) random variable $y$, defined on a common probability space, is defined as $I(x;y)=H(x)-H(x|y)\,,$ where $H(x)$ is the entropy of $x$ (differential entropy if $x$ is a continuous random variable), and $H(x|y)$ is the conditional entropy of $x$ given $y$ ($h(x|y)$ is the conditional differential entropy if $x$ is a continuous random variable). For more general settings including when the random variables are continuous, discrete or a mixture of the two, mutual information is defined as $I(x;y):= \sup_{Q_1,Q_2} I(Q_1(x);Q_2(y)),$ where $Q_1$ and $Q_2$ are quantizers with finitely many bins (see Chapter 5 in \cite{GrayInfo}). An important relevant result is the following. Let $x$ be a random variable and $Q$ be a quantizer applied to $x$. Then, $H(Q(x)) = I(x;Q(x)) = h(x) - h(x | Q(x))$. For a concise overview of relevant information theoretic concepts, we refer the reader to {\sl Chapter 5} of \cite{YukselBasarBook}. For a more complete coverage, see \cite{ElGamalKim} or \cite{Cover}. When the realization $x$ of a random variable $x_t$ needs to be explicitly mentioned, the event $x_t=x$ will be emphasized. We use the conditional probability (expectation) notation $P_x(\cdot)$ ($E_x[\cdot]$)to denote $P(\cdot|x_0=x)$ ($E[\cdot|x_0=x]$). Finally, for a square matrix A, $|A|$ denotes the absolute value of its determinant.

Throughout the paper, all the random variables will be defined on a common probability space $(\Omega, {\cal F}, P)$. 

\subsection{Literature review}

In the literature, the study of non-linear systems have typically considered noise-free controlled systems controlled over discrete noiseless channels. Many of the studies on control of non-linear systems over communication channels have focused on constructive schemes (and not on converse theorems), primarily for noise-free sources and channels, see. e.g. \cite{baillieul2004data}, \cite{liberzon2005stabilization}, and \cite{de2004stabilizability}. For noise-free systems, it typically suffices to only consider a sufficiently small {\it invariant} neighborhood of an equilibrium point to obtain stabilizability conditions. 

One important problem which has not yet been addressed to our knowledge is to obtain converse (or impossibility) theorems: The question of when an open-loop unstable non-linear stochastic control system can or cannot be made ergodic or asymptotically mean stationary subject to information constraints has not been addressed. 

Entropy based arguments (which are crucial in obtaining fundamental bounds in information theory and ergodic theory) can be used to obtain converse results: The entropy, as a measure of uncertainty growth, of a dynamical system has two related interpretations: A topological (distribution-free / geometric) one and a measure-theoretic (probabilistic) one. Although the analysis in this paper is probabilistic, we provide a short discussion on the topological entropy: The distribution-free entropy notion (see, e.g. \cite{katok2007fifty}) for a dynamical system taking values in a compact metric space is concerned with the time-normalized number of distinguishable paths/orbits by some finite $\epsilon > 0$ the system's paths can take values in as the time horizon increases and $\epsilon \to 0$. With such a distribution-free setup \cite{nair2004topological} studied the stabilization of deterministic systems controlled over discrete noiseless finite capacity channels: The topological entropy gives a measure of the number of distinct control inputs needed to make a compact set invariant for a noise-free system. \cite{nair2004topological} extends the notion of topological entropy to controlled dynamical systems, and develops the notion of feedback entropy or invariance entropy \cite{colonius2011invariance}, see also \cite{colonius2009invariance} for related results. \cite{nair2004topological} defines two notions of invariance for a set $K$. A set can be made {\it weakly invariant} if there exists $t > 0$, such that for every $x_0 \in K$, there exists a sequence of control actions so that $x_t \in K' \subset \mbox{interior}(K)$. Strong invariance of $K$ requires that $x_1 \in K'$.  With a relaxation of deterministic controls, \cite{da2013invariance} has studied {\it invariance entropy} for random dynamical systems, and \cite{MatveevSavkin} has generalized the topological entropy theoretic results to include random dynamical models to obtain an observability condition over discrete channels. For a comprehensive discussion of such a geometric interpretation of entropy in controlled systems, see \cite{kawan2013invariance}. The results for deterministic systems pose questions on set stability which are not sufficient to study stochastic setups. Stochasticity also allows for control over general noisy channels, and thus applicable to establish connections with information theory (we note that a distribution-free counterpart for such studies requires one to investigate zero-error capacity formulations \cite{MatveevSavkin}, however many practical channels including erasure channels, have zero zero-error capacity). 

On the other hand, the measure-theoretic (also known as Kolmogorov -- Sinai or metric) entropy is more relevant to information-theoretic as well as random noise-driven stochastic contexts since in this case, one considers the {\it typical} distinguishable paths/orbits of a dynamical system and not all of the sample paths a dynamical system may take (and hence the topological entropy typically provides upper bounds on the measure-theoretic entropy). Measure-theoretic entropy is crucial in the celebrated Shannon-McMillan-Breiman theorem \cite{GrayProbabilit} as well as the isomorphism theorem \cite{Ornstein} \cite{katok2007fifty}. For further relations between different interpretations of entropy as well as their computations (such as through {\it Lyapunov exponents} as a result of Pesin's formula), we refer the reader to \cite{young2003entropy}. Such an entropy notion has operational practical usage in identifying fundamental limits on source and channel coding for stationary sources \cite{ShieldsIT}. However, the findings in the information theory literature has not yet been successfully applied to non-linear networked control systems in general due to the following reasons: (i) The open-loop system in networked control may be unstable and stabilizable only through a control loop. In the information theory literature, stochastic stability results for coding schemes have been established primarily for (control-free) stable sources and when non-stationary, have involved only linear Gaussian auto-regressive (AR) processes \cite{GrayHashimoto}. Moreover, such a control-free analysis does not lead to conclusive results for non-linear controlled sources since non-linear systems suffer from the {\it dual-effect}: one cannot decouple estimation from control, and control from conditional entropy properties under a stationary probability measure. (ii) The coding schemes for such studies in information theory are non-causal; in networked control systems, coding must be causal (that is, real-time or essentially zero-delay \cite{YukselBasarBook}).  

There have been few studies which have adopted a measure-theoretic entropic view for the control of non-linear dynamical systems over communication channels. Relevant contributions include \cite{Mehta1} and \cite{Mehta2}: Building on \cite{Martins2} and \cite{zang2003nonlinear}; \cite{Mehta1}  develops an entropy analysis for non-linear system dynamics to obtain the relation between the entropy rates of a measurement disturbance, output and the dynamical system, and generalizing a {\it Bode-type} entropy analysis for non-linear systems. A related entropy analysis for a class of stochastic non-linear systems have been considered in \cite{Mehta2}. Recently \cite{VaidyaElia} and \cite{vaidya2012stabilization} have considered fading and erasure channels between the controller and the actuator and have studied ergodicity properties using Lyapunov theoretic arguments under a class of structures imposed on control policies; these contributions do not consider finite-rate information and coding restrictions which may arise due to the presence of a channel. Other important relevant work which consider deterministic systems are \cite{liberzon2005stabilization} and \cite{liberzon2002stabilization}, where stability of zooming schemes, as in \cite{BrockettLiberzon}, have been considered. 

Finally, we note an important related discussion in view of Bode's integral formula as extended to a class of non-linear systems in \cite{zang2003nonlinear} under somewhat restrictive conditions, see \cite[Thm.~9]{zang2003nonlinear}. Relevant work includes \cite{EliaTAC}, \cite{Martins2} and \cite{okano2009characterization} for linear systems. For non-linear systems the entropy and mutual information arguments provide the appropriate fundamental bounds instead of a sensitivity integral/transfer function analysis which is commonly used for linear systems as is also advocated in \cite{zang2003nonlinear}. An earlier contribution utilizing measure theoretic entropy for the study and classification of controlled stochastic systems is \cite{Jonckheere}. The findings in our paper provide further generalizations; see Theorem \ref{BodeTheorem} and Remark \ref{BodeRemark}. 

The stability criteria outlined earlier have been studied extensively for linear systems of the form
\begin{eqnarray}\label{vectorEqn1}
x_{t+1}=Ax_t + Bu_t + G w_t,
\end{eqnarray}
where $x_t \in \mathbb{R}^N$ is the state at time $t$, $u_t \in \mathbb{R}^m$ is the control input, and $\{w_t\}$ is a sequence of i.i.d. $\mathbb{R}^{d}$-valued random vectors (such as Gaussian). Here, $(A,B)$ and $(A,G)$ are controllable pairs. %Assume that all eigenvalues $\{\lambda_i, 1 \leq i \leq N\}$ of $A$ are unstable, that is have magnitudes greater than or equal to 1. %There is no loss here since if some eigenvalues are stable, by a similarity transformation, the unstable modes can be decoupled from the stable ones and one can instead consider a lower dimensional system; stable modes are already stochastically stable.

For noise-free linear systems controlled over discrete-noiseless channels, Wong and Brockett \cite{Brockett}, Baillieul \cite{Baillieuil99}; and more generally, Tatikonda and Mitter \cite{Tatikonda} (see also \cite{TatikondaThesis}) and Nair and Evans \cite{NairEvans} have obtained the minimum lower bound needed for stabilization over a class communication channels under various assumptions on the system noise and channels; sometimes referred to as a {\em data-rate theorem}.  This theorem states that for stabilizability under information constraints, in the mean-square sense, a minimum average rate per time stage needed for stabilizability has to be at least $\sum_{i: |\lambda_i| > 1} \log_2(|\lambda_i|)$, where $\{\lambda_i, 1 \leq i \leq N\}$ are the eigenvalues of $A$. 

The particular notion of stochastic stability is crucial in characterizing the conditions on the channels and important extensions have been made in the literature notably by Matveev and Savkin \cite{SavkinSIAM09} \cite{MatveevSavkin}, Sahai and Mitter \cite{Sahai} \cite{SahaiParts}, and Martins et al. \cite{Martins}. For a more comprehensive review; see \cite{NairFagnaniSurvey}, {\sl Chapters 5-8} of \cite{YukselBasarBook}, \cite{Martins2}, and \cite{franceschetti2014elements}. Reference \cite{Minero} considered erasure channels and obtained necessary and sufficient time-varying rate conditions for control over such channels. Reference \cite{Coviello} considered second moment stability over a class of Markov channels with feedback. Motivated from such problems, \cite{YukTAC2010} and \cite{YukMeynTAC2010} developed a martingale-method for establishing stochastic stability, which later led to a  random-time state-dependent drift criterion, leading to the existence of an invariant distribution possibly with moment constraints; these were utilized to obtain policies leading to strong forms of stochastic stability, such as ergodicity or positive Harris recurrence \cite{YukselBasarBook}, for linear systems driven by additive unbounded noise.

The following definition (see \cite[Definition 8.5.1]{YukselBasarBook}) will be useful in the analysis later in the paper.
\begin{defn}\label{ClassAChannelDefinition}
Channels are said to be of Class A type, if
\begin{itemize}
\item they satisfy the following Markov chain condition:
\begin{eqnarray}
 q'_t \leftrightarrow q_t, q_{[0,t-1]}, q'_{[0,t-1]} \leftrightarrow \{x_{0}, w_s, s \geq 0\},\label{classAMarkov}
 \end{eqnarray}
 that is, almost surely, for all Borel sets $B$, 
 \[P(q'_t \in B| q_t, q_{[0,t-1]}, q'_{[0,t-1]}, x_{0}, w_s, s \geq 0) = P(q'_t \in B| q_t, q_{[0,t-1]}, q'_{[0,t-1]}),\]
 for all $t \geq 0$, and
\item their capacity with feedback is given by:
\[C= \lim_{T \to \infty} \max_{\{ P(q_t | q_{[0,t-1]},q'_{[0,t-1]}), \quad 0 \leq t \leq T-1\}} {1 \over T} I(q_{[0,T-1]} \to q'_{[0,T-1]}),\]
where the directed mutual information is defined by \[I(q_{[0,T-1]} \to q'_{[0,T-1]}) = \sum_{t=1}^{T-1} I(q_{[0,t]};q'_t|q'_{[0,t-1]}) + I(q_0;q'_0).\]
%with conditioning on $q'_{-1}$ is meant to be conditioning on the (trivial) $\sigma-field$ generated by $\emptyset$.
\end{itemize}
\end{defn}

Memoryless channels belong to this class; for such channels, feedback does not increase the capacity \cite{Cover}. Such a class also includes finite state stationary Markov channels which are indecomposable \cite{PermuterWeissmanGoldsmith}, and non-Markov channels which satisfy certain symmetry properties \cite{SenAlaYukIT}. Further examples can be found in \cite{TatikondaIT} and in \cite{DaboraGoldsmith}.

\begin{thm}\label{StochasticStabilityofVectorNecessity} \cite{YukselBasarBook}  \cite{YukselAMSITArxiv}
Consider the multi-dimensional linear system (\ref{vectorEqn1}). For such a system controlled over a Class A type noisy channel with feedback, if the channel capacity satisfies
\[C < \sum_{|\lambda_i|>1}  \log_2(|\lambda_i|),\]
(i) there does not exist a stabilizing coding and control scheme with the property $\liminf_{T \to \infty} {1 \over T} {h(x_T)} \leq 0$, (ii) the system cannot be made AMS or ergodic (see Section \ref{ergodicTheory}).
 \end{thm}

For sufficiency, assume that $A$ is a diagonalizable matrix (a sufficient condition for which is that its eigenvalues are distinct real).

\begin{thm}\label{StochasticStabilityofVector} \cite{YukselBasarBook} \cite{YukselAMSITArxiv}
Consider the multi-dimensional linear system (\ref{vectorEqn1}) with a diagonalizable matrix $A$ and Gaussian noise, controlled over a discrete memoryless channel. If the Shannon capacity of the channel satisfies
\[C > \sum_{|\lambda_i|>1} \log_2(|\lambda_i|),\]
there exists a stabilizing scheme which makes the process $\{x_t\}$ AMS. If the channel is noiseless, or a memoryless erasure channel, or is a Gaussian channel, then the process $\{x_t\}$ can be made stationary and ergodic.
\end{thm}

\subsection{Contributions of the paper}

As stated above, stochastic stabilization of non-linear systems driven by noise (especially unbounded noise) over communication channels has not been studied to our knowledge where the goal is to establish asymptotic (mean) stationarity, ergodicity or stationarity of the closed-loop system. We use measure-theoretic entropy analysis and ergodic theoretic tools for arrive at necessary and sufficient conditions. A by-product of the analysis is a generalization of Bode's Integral Formula to a class of non-linear systems and arbitrary information channels with memory. The approach in the paper, although building on our earlier work on linear systems, contains significant generalizations in the approach due to the non-linearity of the source. We also consider a construction of a stabilizing coding and control scheme for multi-dimensional non-linear sources driven by unbounded noise controlled over a discrete noiseless channel.

\section{Sublinear entropy growth and a generalization of Bode's Integral Formula for non-linear systems}\label{Bode}

In the paper, instead of a general $\mathbb{R}^N$-valued non-linear state model
\begin{eqnarray}
x_{n+1} = f(x_n,u_n,w_n),\label{model4}
\end{eqnarray}
we will consider non-linear systems of the form
\begin{eqnarray}
x_{n+1} &=& f(x_n, w_n) + B u_n,\label{model2} \\
x_{n+1} &=& f(x_n) + Bu_n + w_n,\label{model3} \\
x_{n+1} &=& f(x_n,u_n) + w_n.\label{model5} 
\end{eqnarray}
We also will have an occasion discuss non-linear systems of the form
\begin{eqnarray}
x_{n+1} = f(x_n, w_n) + B(x_n) u_n. \label{model1}
\end{eqnarray}

In all of the models above, $x_n$ is the $\mathbb{R}^{N}$-valued state, $w_n$ is the $\mathbb{R}^N$-valued noise variable, $u_n$ is $\mathbb{R}^s$ valued and $w_n$ assumed to be an independent noise process with $w_n \sim \nu$. %Such models include stochastic control-affine systems as well as bi-linear control systems. 

We assume throughout that $f$ is measurable and continuously differentiable in the state variable. For a possibly non-linear differentiable function $f: \mathbb{R}^N \to \mathbb{R}^m$, the Jacobian matrix of $f$ is an $n \times m$ matrix function consisting of partial derivatives of $f$ such that \[J(f)(i,j) = { \partial (f(x))_i \over \partial x_j}, \quad 1 \leq i \leq m, 1 \leq j \leq n.\] 

We will have the following assumption throughout the paper.
\begin{assumption}\label{invertiblef}
In the models considered above $f(\cdot,w): \mathbb{R}^N \to \mathbb{R}^N$ is invertible for every realization of $w$.
\end{assumption}

In the following $|J(f)|$ will denote the absolute value of the determinant of the Jacobian. Furthermore, with $f_w(x) = f(x,w)$, we define $J(f(x,w)) := J(f_w(x))$. % instead of $J(f(x,w))$ to highlight the differentiation in the $x$ variable only.

\begin{assumption}\label{bounds2}
There exist $M_1 \in \mathbb{R}$ and $L_1 \in \mathbb{R}$ so that for all $x,w$
\[L_1 \leq \log_2(|J(f(x,w))|) \leq M_1\]
\end{assumption}

The following is our first result; it provides conditions for sublinear entropy growth (in time) which implies quadratic stability. The result will also be used in the next section and its proof leads to a generalization of Bode's Integral Formula as discussed further below. Let $\pi_t(B) = P(x_t \in B)$ for all Borel $B$.

\begin{thm}\label{sublinear}
Consider the networked control problem over a Class A channel.
(i) Let f have the form in (\ref{model2}), (ii) Assumptions \ref{invertiblef} and \ref{bounds2} hold, and (iii) $x_0$ have finite differential entropy. 
a) If there is an admissible coding and control policy such that
% \[\liminf_{t \to \infty} h(x_t|q'_{[0,t-1]})/t \leq 0,\]
%or
\[\liminf_{t \to \infty} h(x_t)/t \leq 0,\]
it must be that 
%$C \geq  \lim_{t \to \infty} h(f(x_t,u_t,w_t) | q'_{[0,t-1]}) - h(f(x_t,u_t,w_t) | q'_{[0,t-1]})$ for some policy $u_t = \gamma_t(q'_{[0,t-1]})$
\begin{eqnarray}
C \geq  \liminf_{T \to \infty} {1 \over T} \sum_{t=0}^{T-1} \int \pi_t(dx) \bigg(\int \nu(dw) \log_2(|J(f(x,w))|) \bigg)  \label{convAvg}
\end{eqnarray}
%where the expectation is over the channel output variables $\{q'_t\}$.
%\[C \geq   \limsup_{T \to \infty} {1 \over T}\bigg( \bigg(\sum_{t=1}^{T-1} \int P(w_{t-1}) \log_2(|J(f(x_{t-1},w_{t-1}))|) \bigg)\]
b)  If there is an admissible coding and control policy such that 
%\[\limsup_{t \to \infty} h(x_t|q'_{[0,t-1]})/t \leq 0,\]
%or
\[\limsup_{t \to \infty} h(x_t)/t \leq 0,\]
it must be that 
%$C \geq  \lim_{t \to \infty} h(f(x_t,u_t,w_t) | q'_{[0,t-1]}) - h(f(x_t,u_t,w_t) | q'_{[0,t-1]})$ for some policy $u_t = \gamma_t(q'_{[0,t-1]})$
\begin{eqnarray}
C \geq  \limsup_{T \to \infty} {1 \over T}\sum_{t=0}^{T-1} \int \pi_t(dx) \bigg(\int \nu(dw) \log_2(|J(f(x,w))|) \bigg)  \label{convAvg2}
\end{eqnarray}
In either case, if $L := \inf_{x,w} \log_2|J(f(x,w))|$, then $C \geq L$.
\end{thm}

\begin{remark}
The condition $\limsup_{t \to \infty} h(x_t)/t \leq 0$ is a weak condition. For example a stochastic process whose second moment grows subexponentially in time so that $\limsup_{T \to \infty} {\log(E[x_T^2]) \over T} \leq 0$, satisfies this condition. Hence, {\it quadratic stability} implies this condition. \hfill $\diamond$
\end{remark}

\begin{remark}
In the theorem, we would have obtained the same results if we had replaced $\limsup_{t \to \infty} h(x_t)/t \leq 0$ with $\limsup_{t \to \infty} {1 \over t}h(x_t|q'_{[0,t-1]})/t \leq 0$. This condition would be more relevant for state
estimation problems, where the goal is not necessarily to make the state stable, but to make the estimation error stable (where $u_t$ would be the state estimate and $x_t - u_t$ would be the estimation error). Since $h(x_t|q'_{[0,t-1]}) \leq h(x_t)$, it is evident that the condition $h(x_t)/t \leq 0$ implies that $h(x_t|q'_{[0,t-1]})/t \leq 0$. \hfill $\diamond$ 
\end{remark}

\textbf{Proof of Theorem \ref{sublinear}}\label{proofNecessity}
Recall that for channels of the type Class A (which includes the discrete memoryless channels (DMC) as a special case), the capacity is given by:
 \[C= \lim_{T \to \infty} \max_{\{P(q_t | q_{[0,t-1]},q'_{[0,t-1]}) \}} {1 \over T} I(q_{[0,T-1]} \to q'_{[0,T-1]})\]
where
\[I(q_{[0,T-1]} \to q'_{[0,T-1]}) = \sum_{t=1}^{T-1} I(q_{[0,t]};q'_t|q'_{[0,t-1]}) + I(x_0; q'_0).\]

Let us define $R_T=  \max_{\{ P(q_t | q_{[0,t-1]},q'_{[0,t-1]}), 0 \leq t \leq T-1\}} {1 \over T}\sum_{t=0}^{T-1} I(q'_t; q_{[0,t]} | q'_{[0,t-1]})$. Observe that for $t > 0$:
\begin{eqnarray}\label{usefulStep}
I(q'_t; q_{[0,t]} | q'_{[0,t-1]}) && = H(q'_t | q'_{[0,t-1]}) -  H(q'_t | q_{[0,t]} ,q'_{[0,t-1]}) \nonumber \\
&& = H(q'_t | q'_{[0,t-1]}) -  H(q'_t | q_{[0,t]}, x_t,q'_{[0,t-1]}) \label{classAassumption} \\
&& \geq H(q'_t | q'_{[0,t-1]}) -  H(q'_t | x_t,q'_{[0,t-1]}) \nonumber \\
&& = I(x_t; q'_t| q'_{[0,t-1]}).
\end{eqnarray}
Here, (\ref{classAassumption}) follows from the assumption that the channel is of Class A type. 
\[\quad\]
a) Consider the following
\begin{eqnarray}
&&\lim_{T \to \infty} R_T \geq \limsup_{T \to \infty} {1 \over T} \bigg( \sum_{t=1}^{T-1} I(x_t; q'_t| q'_{[0,t-1]})) + I(x_0;q'_0) \bigg) \label{usefulStep2} \\
%&\geq& \limsup_{T \to \infty} {1 \over T}\bigg( \sum_{t=1}^{T-1} \bigg(H(q'_t | q'_{[0,t-1]}) - H(q'_t |x_t, q'_{[0,t-1]})\bigg) + I(x_0;q'_0) \bigg)\nonumber \\
%&=& \limsup_{T \to \infty} {1 \over T}\bigg( \sum_{t=1}^{T-1} \bigg(I(x_t;q'_t | q'_{[0,t-1]})\bigg) +I(x_0;q'_0)\bigg)\nonumber \\
 &=&  \limsup_{T \to \infty} {1 \over T}\bigg( \sum_{t=1}^{T-1} \bigg(h(x_t | q'_{[0,t-1]}) - h(x_t|q'_{[0,t]})\bigg) + I(x_0;q'_0)\bigg)\nonumber \\
 &=&  \limsup_{T \to \infty} {1 \over T}  \sum_{t=1}^{T-1} \bigg(h(x_t | q'_{[0,t-1]}) - h(x_t|q'_{[0,t]})\bigg) \nonumber \\
 &=&  \limsup_{T \to \infty} {1 \over T}  \sum_{t=1}^{T-1} \bigg(h(f(x_{t-1},w_{t-1}) + B u_{t-1} | q'_{[0,t-1]}) - h(x_t|q'_{[0,t]})\bigg) \nonumber \\
 &=&  \limsup_{T \to \infty} {1 \over T}  \sum_{t=1}^{T-1} \bigg(h(f(x_{t-1},w_{t-1}) | q'_{[0,t-1]}) - h(x_t|q'_{[0,t]})\bigg) \nonumber \\
 &=&  \limsup_{T \to \infty} {1 \over T}  \sum_{t=1}^{T-1} \bigg( \sum_{\zeta_{[0,t-1]}} h(f(x_{t-1},w_{t-1}) | q'_{[0,t-1]} = \zeta_{[0,t-1]}) P(q'_{[0,t-1]} = \zeta_{[0,t-1]})   \nonumber \\
&& \quad \quad \quad \quad \quad \quad \quad \quad- h(x_t|q'_{[0,t]})\bigg) \label{DefnCondition1} \\
&\geq&  \limsup_{T \to \infty} {1 \over T}  \sum_{t=1}^{T-1} \bigg( \sum_{\zeta_{[0,t-1]}} h(f(x_{t-1},w_{t-1}) | q'_{[0,t-1]} = \zeta_{[0,t-1]},w_{t-1})  P(q'_{[0,t-1]} = \zeta_{[0,t-1]}) \bigg) \nonumber \\ && \quad \quad \quad \quad \quad \quad \quad  \quad \quad  - h(x_t|q'_{[0,t]})  \label{conditionW} \\
&=&  \limsup_{T \to \infty} {1 \over T} \sum_{t=1}^{T-1} \bigg( \bigg( \sum_{\zeta_{[0,t-1]}} \int h(f(x_{t-1},w) | q'_{[0,t-1]} = \zeta_{[0,t-1]},w_{t-1}=w) v(dw) \nonumber \\ && \quad \quad \quad \quad \quad \quad \quad \quad \quad  \quad \quad \times P(q'_{[0,t-1]} = \zeta_{[0,t-1]}) \bigg)     - h(x_t|q'_{[0,t]})\bigg)  \label{DefnCondition2} \\
 &=&  \limsup_{T \to \infty} {1 \over T}  \sum_{t=1}^{T-1}  \bigg( \sum_{\zeta_{[0,t-1]}} P(q'_{[0,t-1]} = \zeta_{[0,t-1]}) \nonumber \\ &&  \quad  \quad \times \bigg( \int \nu(dw) \bigg( \int P(dx_{t-1}|q'_{[0,t-1]}= \zeta_{[0,t-1]},w_{t-1}=w) \log_2(|J(f(x_{t-1},w))|) \nonumber \\ &&  \quad \quad \quad \quad  \quad \quad \quad \quad \quad  + h(x_{t-1} | q'_{[0,t-1]}= \zeta_{[0,t-1]},w_{t-1}=w) \bigg) \bigg) - h(x_t|q'_{[0,t]})\bigg) \label{indepW0} \\
%&& \quad \quad \quad \quad \quad  - h(x_t|q'_{[0,t]})\bigg) + I(x_0;q'_0)\bigg)  \label{indepW0} \\
 &=&  \limsup_{T \to \infty} {1 \over T}  \sum_{t=1}^{T-1}  \bigg( \sum_{\zeta_{[0,t-1]}} P(q'_{[0,t-1]} = \zeta_{[0,t-1]}) \nonumber \\ &&  \quad \quad \quad \quad  \quad \quad \times \bigg( \int \nu(dw) \bigg( \int P(dx_{t-1}|q'_{[0,t-1]}= \zeta_{[0,t-1]})  \log_2(|J(f(x_{t-1},w))|)  \nonumber \\
&& \quad \quad \quad \quad \quad + h(x_{t-1} | q'_{[0,t-1]}= \zeta_{[0,t-1]}) \bigg)  \bigg) - h(x_t|q'_{[0,t]})\bigg) \label{indepWv} \\ 
 &=&  \limsup_{T \to \infty} {1 \over T}   \sum_{t=1}^{T-1} \bigg( \bigg( \sum_{\zeta_{[0,t-1]}} P(q'_{[0,t-1]} = \zeta_{[0,t-1]}) \int P(dx_{t-1}|q'_{[0,t-1]} = \zeta_{[0,t-1]} ) \nonumber \\ &&   \quad  \quad \quad \times  \int \nu(dw) \log_2(|J(f(x_{t-1},w))|) \bigg) + h(x_{t-1} | q'_{[0,t-1]}) - h(x_t|q'_{[0,t]})\bigg) \label{indepWv2} \\
 &=&  \limsup_{T \to \infty} {1 \over T} \bigg( \sum_{t=0}^{T-1} \int \pi_t(dx) \bigg(\int \nu(dw) \log_2(|J(f(x,w))|) \bigg) - h(x_{T-1} | q'_{[0,T-1]}) \bigg) \nonumber \\
% &=&  \limsup_{T \to \infty} {1 \over T}\bigg( \sum_{t=1}^{T-1} \bigg(\int_{w_{t-1}} P(dw_{t-1}) \bigg( \int P(dx_{t-1}|q'_{[0,t-1]}) \log_2(|J(f(x_{t-1},w_{t-1})|) + h(x_{t-1} | q'_{[0,t-1]}) - h(x_t|q'_{[0,t]})\bigg) + I(x_0;q'_0)\bigg)\nonumber \\
% &=&  \limsup_{T \to \infty} {1 \over T}\bigg( \sum_{t=1}^{T-1} \bigg(\int_{w_{t-1}} P(dw_{t-1}) \bigg(\int P(dx_{t-1}|q'_{[0,t-1]})  \log_2(|J(f(x_{t-1},w_{t-1})|) + h(x_{t-1} | q'_{[0,t-1]}) - h(x_t|q'_{[0,t]})\bigg) + I(x_0;q'_0)\bigg)\nonumber \\
 &\geq&  V - \liminf_{T \to \infty}  {1 \over T} h(x_{T-1} | q'_{[0,T-1]})  \label{boundeNoAMS1} 
% &=&  \log_2(|A|) - \liminf_{T \to \infty} \bigg( {1 \over T} h(x_{T-1} - \bar{x}_t(q'_{[0,T-1]}) | q'_{[0,T-1]}) \bigg)  \nonumber
%&\geq&  V - \liminf_{T \to \infty} \bigg( {1 \over T} h(x_{T-1}) \bigg) \geq V. \label{bounde} \label{conditionagain}
\end{eqnarray}
Here, 
 \begin{eqnarray}
%V = \liminf_{T \to \infty} {1 \over T}E\bigg[\bigg( \sum_{t=1}^{T-1} \int \nu(dw) \int P(dx_{t-1}|q'_{[0,t-1]}) \log_2(|J(f(x_{t-1},w))|) \bigg) \bigg] \nonumber \\ \label{Vdefn}
V := \liminf_{T \to \infty} {1 \over T} \sum_{t=0}^{T-1} \int \pi_t(dx) \bigg(\int \nu(dw) \log_2(|J(f(x,w))|) \bigg) \label{Vdefn}
%V= \liminf_{T \to \infty} {1 \over T} \bigg( \sum_{t=1}^{T-1} \int \nu(dw) \int \pi_{t-1}(dx_{t-1}) \log_2(|J(f(x_{t-1},w))|) \bigg) \label{Vdefn}
 \end{eqnarray}
% where we take out the conditioning on $q'_{[0,t-1]}$ using the properties of iterated expectations.

Equations (\ref{DefnCondition1}) and (\ref{DefnCondition2}) follow from the definition of conditional entropy, (\ref{conditionW}) follows from conditioning on the random variable $w_{t-1}$. Equations (\ref{indepW0})-(\ref{indepWv}) follow from the fact that $x_{t-1} \leftrightarrow q'_{[0,t-1]} \leftrightarrow w_{t-1}$ is a Markov chain and the following. For every realization $q'_{[0,t-1]} = \zeta_{[0,t-1]}$,
\begin{eqnarray} 
&& h(f(x_{t-1},w_{t-1}) | q'_{[0,t-1]} = \zeta_{[0,t-1]},w_{t-1}=w) \nonumber \\
&& = h(f^w(x_{t-1}) | q'_{[0,t-1]} = \zeta_{[0,t-1]},w_{t-1}=w) \nonumber \\
&& = \int P(dx_{t-1}|q'_{[0,t-1]} = \zeta_{[0,t-1]},w_{t-1}=w) \log_2(|J(f_w(x_{t-1}))|) \nonumber \\
&& \quad \quad \quad \quad + h(x_{t-1} | q'_{[0,t-1]} = \zeta_{[0,t-1]},w_{t-1}=w) \nonumber \\
 \label{entropyTranslate} \\
&& = \int P(dx_{t-1}|q'_{[0,t-1]} = \zeta_{[0,t-1]}) \log_2(|J(f(x_{t-1},w))|) + h(x_{t-1} | q'_{[0,t-1]} = \zeta_{[0,t-1]}), \nonumber
\end{eqnarray} 
where $f_w(x) := f(x,w)$ is an invertible function for every $w$, and as a result (\ref{entropyTranslate}) follows from the entropy formula for invertible functions of a random variables (see, e.g., p. 167 of \cite{StarkWoods} and Lemma 4 in \cite{zang2003nonlinear}) and the last line follows from the condition $x_{t-1} \leftrightarrow q'_{[0,t-1]} \leftrightarrow w_{t-1}$. Equation (\ref{indepWv2}) follows from Fubini's theorem by Assumption \ref{bounds2}.

By the hypothesis,
 $\liminf_{t \to \infty} {1 \over t} h(x_t) \leq 0$,
 it must be that $\lim_{T \to \infty} R_T \geq V$. Thus, the capacity also needs to satisfy this bound. 
 
 In the above derivation, (\ref{boundeNoAMS1}) follows from the fact that for two sequences $a_n, b_n$:
\begin{eqnarray}\label{limsupliminfineq}
\limsup_{n \to \infty} (a_n + b_n) \geq \limsup_{n \to \infty} a_n + \liminf_{n \to \infty} b_n.
\end{eqnarray}

b) If $\limsup_{t \to \infty} h(x_t|q'_{[0,t-1]})/t \leq 0,$ (\ref{boundeNoAMS1}) can be applied through (\ref{limsupliminfineq}) with $V$ defined as
\[\limsup_{T \to \infty} {1 \over T}E\bigg[ \bigg( \sum_{t=1}^{T-1} \int P(dx_t| q'_{[0,t-1]})  \bigg(\int \nu(dw) \log_2(|J(f(x_t,w))|) \bigg) \bigg) \bigg] \]
and in (\ref{boundeNoAMS1}), $\liminf_{T \to \infty} {1 \over T} h(x_{T-1} | q'_{[0,T-1]})$ being replaced with $\limsup$ of the same expression.
%_{T \to \infty} {1 \over T} h(x_{T-1} | q'_{[0,T-1]})$.
\qed

\begin{remark}
We note that if the system had been of a model in (\ref{model1}), the expression involving $J(f(x,w))$ would explicitly depend on the control policy which would in turn depend possibly on the entire past channel outputs making the expression computationally more involved. \hfill $\diamond$%We note also that $P(dx_t|q'_{[0,t-1]})$ depends implicitly on the control policy adopted and despite this fact, the theorem provides explicit quantitative bounds which will be discussed further later. \hfill $\diamond$
\end{remark}

\subsection{A generalization of Bode's Integral Formula for non-linear systems}
The proof of Theorem \ref{sublinear} reveals an interesting connection with and generalization of Bode's Integral Formula (and what is known as the {\it waterbed effect}) \cite{middleton1991trade} to non-linear systems, which we state formally in the following. The result also suggests that an appropriate generalization for non-linear systems is through an information theoretic approach that recovers Bode's original result for the linear case as we discuss further below. 
\begin{thm}\label{BodeTheorem}
(i) Let f have the form in (\ref{model2}), (ii) Assumption \ref{invertiblef} hold, and (iii) $x_0$ have finite differential entropy. 
If there is an admissible coding and control policy with $\limsup_{t \to \infty} h(x_t)/t \leq 0$
% \[\limsup_{t \to \infty} h(x_t|q'_{[0,t-1]})/t \leq 0, \quad $\mathrm{or}$ \quad \limsup_{t \to \infty} h(x_t)/t \leq 0,\]
it must be that 
\begin{eqnarray}
&& \limsup_{T \to \infty} {1 \over T} I(q_{[0,T-1]} \to q'_{[0,T-1]}) \nonumber \\
&& \quad \quad \geq \limsup_{T \to \infty} {1 \over T}  \sum_{t=0}^{T-1} \int \pi_t(dx) \bigg(\int \nu(dw) \log_2(|J(f(x,w))|) \bigg) \label{convAvg2}
%E\bigg[ \bigg( \sum_{t=1}^{T-1} \int P(x_t \in dx| q'_{[0,t-1]}) \int \nu(dw) \log_2(|J(f(x,w))|) \bigg) \bigg] \nonumber \\ \label{convAvg2}
\end{eqnarray}
%\[C \geq   \limsup_{T \to \infty} {1 \over T}\bigg( \bigg(\sum_{t=1}^{T-1} \int P(w_{t-1}) \log_2(|J(f(x_{t-1},w_{t-1}))|) \bigg)\]
%In particular if $L := \inf_{x,w} \log_2|J(f(x,w))|$, then $C \geq L$.
\end{thm}
\textbf{Proof.} This follows directly from equations (\ref{usefulStep}),(\ref{usefulStep2}) and (\ref{boundeNoAMS1}).  \qed

\begin{remark}\label{BodeRemark} {\bf [Reduction to Bode's Integral Formula for Linear Systems and Gaussian Noise]} If the system considered is linear with all open-loop eigenvalues unstable, the channel is an additive noise channel so that $q'_t=q_t+v_t$ for some stationary Gaussian noise, and time-invariant control policies are considered leading to a stable system, then with the more common notation of $y_t=q'_t$, the right hand side of (\ref{convAvg2}) would be the sum of the unstable eigenvalues of the linear system matrix. For a stationary Gaussian process [see \cite{Cover}, page 274] the entropy rate can be written as \[{1 \over 2} \log(2 \pi e) + \int_{-1/2}^{1/2} {1\over 2} \log(S(f)) df\]
 with $S$ denoting the spectral density of the process. Now, (\ref{usefulStep}) becomes
\[I(q'_t; q_{[0,t]} | q'_{[0,t-1]})  = h(q'_t | q'_{[0,t-1]}) -  h(q'_t | q_{[0,t]} ,q'_{[0,t-1]}) = h(q'_t | q'_{[0,t-1]}) -  h(v_t | v_{[0,t-1]}),\]
and thus the left hand side of (\ref{convAvg2}) reduces to the difference between the entropy rate of the process $q'_t$ (that is, $\lim_{t \to \infty} h(q'_t | q_{[0,t-1]})$) and that of the stationary noise process $v_t$  (that is, $\lim_{t \to \infty} h(v_t | v_{[0,t-1]})$. Then, the left hand side of (\ref{convAvg2}) equals 
 \[\int_{-1/2}^{1/2} {1\over 2} \log({S_{y}(f) \over S_v(f)}) df,\] 
 which then is equal to the integral of the log-sensitivity function (corresponding to the transfer function from the disturbance process $v_t$ to the output process $q'_t$). 
%Furthermore, (\ref{usefulStep}) becomes
%$I(q'_t; q_{[0,t]} | q'_{[0,t-1]})  = H(q'_t | q'_{[0,t-1]}) -  H(q'_t | q_{[0,t]} ,q'_{[0,t-1]}) = H(q'_t | q'_{[0,t-1]}) -  H(v_t | v_{[0,t-1]})$ and follows that the inequalities in the proof of Theorem \ref{sublinear} become inequalities and we obtain the standard form of Bode's formula:
%\[\int_{-1/2}^{1/2} {1\over 2} \log({S_{y}(f) \over S_v(f)}) df = \int \log_2(|\lambda_i|).\] 
 This leads to the celebrated Bode's Integral Formula. In the context of linear systems, earlier extensions of this formula have been studied in \cite{EliaTAC} with an information theoretic interpretation under the restriction to linear policies  (see e.g. Theorem 4.6 in \cite{EliaTAC}), in \cite{Martins2} under more general possibly non-linear stabilizing control policies which lead to a stationary process, and in \cite{Mehta1} and \cite{zang2003nonlinear} for a class of non-linear noise-free systems. \hfill $\diamond$
\end{remark}

\section{Asymptotic mean stationarity and ergodicity}\label{Ergodicity}

In the following, we build on, but significantly modify the approaches in \cite{Matveev} and \cite{YukselBasarBook} to account for non-linearity of the system.

Consider the system (\ref{model3}), under some admissible policy, controlled over a channel. 
\begin{assumption}\label{bounds1}
We assume
\begin{eqnarray}
M &:=& \sup_{x \in \mathbb{R}^N} \log_2|J(f(x))| <  \infty, \nonumber \\
L &:=& \inf_{x \in \mathbb{R}^N} \log_2|J(f(x))| > - \infty. \nonumber
\end{eqnarray}
\end{assumption}

\begin{thm}\label{ModifiedMatveevsLemma} 
Consider the system (\ref{model3}) controlled over a Class A type noisy channel with feedback where $h(x_0) < \infty$ and Assumptions \ref{invertiblef} and \ref{bounds1} hold. If $C < L$, then under any admissible policy,
\[\limsup_{T \to \infty} P(|x_T| \leq b(T)) \leq 1 - {L-C \over M},\]
for all $b(T) > 0$ such that $\lim_{T \to \infty} {1 \over T} \log_2(b(T)) = 0$.
\end{thm}

The proof is in Section \ref{ProofofModified} of the Appendix. An implication of this result follows.
\begin{thm}\label{notAMS}
Consider the system (\ref{model3}) controlled over a Class A type noisy channel with feedback where $h(x_0) < \infty$ and Assumptions \ref{invertiblef} and \ref{bounds1} hold. If under some causal encoding and controller policy the state process is AMS, the channel capacity $C$ must satisfy $C \geq L$.
\end{thm}

We recover the following result for linear systems in \cite{YukselBasarBook} as a special case.
\begin{cor}
For the linear case with $f(x)=Ax$ with eigenvalues $|\lambda_i| \geq 1$, $C \geq \sum_k \log_2(|\lambda_i|)$ is a necessary condition for the AMS property under any admissible coding and control policy.
\end{cor}

\textbf{Proof of Theorem \ref{notAMS}} \label{proofofnotAMSExtension}
If the process is AMS (see Section \ref{ergodicTheory}), then there exists a stationary measure $\bar{P}$ such that
\begin{eqnarray}\label{AMSconverseStep}
\lim_{N \to \infty} {1 \over N} \sum_{k=1}^N P(T^{-k}D) = \bar{P}(D),
\end{eqnarray}
for all (cylinder) events $D$. Let for $b_B \in \mathbb{R}_+$, $B \in {\cal B}(\mathbb{R}^N)$ be given by $B = \{x: |x| \leq b_B\}$ and $X_n(z)=z_n$ be the coordinate function (see Section \ref{ergodicTheory}) where $z = \{z_0, z_1, z_2, \cdots\}$.
%Let given $n$, $B=\{x: |X_n(x)| \leq b_B\}$, where $B = [-b_B,b_B]$ and $X_n(z)=z_n$ be the coordinate function (see Section \ref{ergodicTheory}).

If by Theorem \ref{ModifiedMatveevsLemma}
\begin{eqnarray}\label{MatveevsCondition8}
\limsup_{T \to \infty} P(|x_T| \leq b_B) \leq 1 - {(L-C) \over M} < 1,
\end{eqnarray}
holds for all $b_B \in \mathbb{R}_+$, then $\bar{P}_n(B) < 1 - {(L-C) \over M}$ for all compact $B$, where $\bar{P}_n$ is the marginal probability on the $n$th coordinate defined as
\[\bar{P}_n(B) = \bar{P}\bigg(x: |X_n(x)| \leq b_B\bigg).\]
But then $\bar{P}_n$, as an individual probability measure, must be tight \cite{Billingsley}, therefore, for every $\delta > 0$ there exists $b_B < \infty$ such that $\bar{P}_n(B) \geq 1 - \delta$. But, by (\ref{AMSconverseStep}), this would imply that $\limsup_{t \to \infty} P(T^{-t}B)=\limsup_{t \to \infty} P(|x_t| \in B) \geq 1 - \delta$, leading to a contradiction with (\ref{MatveevsCondition8}) for $\delta <  {L - C \over M}$. Hence, the AMS property cannot be achieved. \qed

We end this section with a remark.
\begin{remark}
In information theory, a well-established result is that for noiseless coding of information stable sources (this includes all finite state stationary and ergodic sources) over a class of information stable noisy channels (which includes the channels we consider here), an asymptotically noise-free recovery is possible if the channel capacity is greater than the source entropy through the use of non-causal codes, see e.g. \cite{VerduSeparation} \cite{kieffer1981zero}. However, for the problem we consider (i) the source is non-stationary and open-loop unstable, (ii) the encoding is causal, and (iii) the source process space is not finite-alphabet. Nonetheless, we see that the invariance properties of the source process does appear in the rate bounds that we obtain. \hfill $\diamond$
\end{remark}
%\begin{remark}
%In topological feedback entropy, the authors either consider stabilization to a set, or a neighborhood around a point. The system is noise-free and it is assumed that for every point inside the set, there exists a finite time $T$ such that, one can find a collection of control actions so that $x_T=f(x_0,u_0,\cdots,u_{T_1})$ is in the set. This condition implies the AMS property; but it is also very strong, the dependency on $T$ is uniform here. For the point stabilization, further assumptions are made on the system dynamics using the Jordan form of the Jacobian matrix.
%\end{remark}

\section{Stationarity and positive Harris recurrence under structured (stationary) policies}\label{Stationary}

In many applications, one uses a state-space formulation for coding and control policies. In the following, we will consider stationary update rules which have the form that 
\begin{eqnarray}
&& q_t = \gamma^e(x_t,m_t) \nonumber \\
&&u_t= \gamma^d(m_t,q'_t), \nonumber \\
&&m_t = \eta(m_{t-1},q'_{t-1}), \label{updates}
\end{eqnarray} 
for functions $\gamma^e, \gamma^d$, and $\eta$. In the form above, $m$ is a $\mathbb{S}$-valued memory or {\it quantizer state} variable. A large class of adaptive encoding policies have this form. This includes, delta modulation, differential pulse coded modulation (DPCM), adaptive differential pulse coded modulation (ADPCM), Goodman-Gersho type adaptive quantizers (see e.g. \cite{Kieffer} \cite{KiefferDunham}), as well as the coding schemes used for stabilization of networked control systems under fixed-rate codes \cite{YukTAC2010}. Even further, jointly  optimal source and channel codes for zero-delay coding schemes under infinite horizon optimization criteria also have the form above (where $\mathbb{S}$ is a space of probability measures \cite{YukLinZeroDelay}). We now present a necessary structural result on the encoders. 
\subsection{A necessary structural result on the encoders}

Let $m_t$, the state of the encoder, take values in $\mathbb{S}$. Consider (\ref{model3}). A stabilizing time-invariant encoder/decoder/controller policy given (\ref{updates}), in general, cannot have $|\mathbb{S}| < \infty$.

\begin{thm}\label{Transience111}
Consider (\ref{model3}) with scalar $x_t$ and $w_t$ with a probability measure $\nu$ such that it has a density positive everywhere and $E_{\nu}[\gamma^{-w}] < \infty$ for some $\gamma > 1$. Suppose that there exists $K > 0$ so that
\[\inf_{x > K} {df \over dx}(x) > 1. \]
and ${df \over dx}(x)$ is bounded. Then, a finite cardinality for $\mathbb{S}$, under (\ref{updates}) leads to a transient system in the sense that
 \[P_x(\tau_S < \infty) < 1\] where for some $s>0$, $S = (-\infty, s)$ is an open set containing the origin, $x > s$ and $\tau_S := \inf(t> 0: x_t \in S)$.
A similar result applies for the condition 
\[\sup_{x < - K} {df \over dx}(x) < -1,\]
with $S = (s,\infty)$ for some $s<0$ and $x < s$.
\end{thm}

\textbf{Proof.}
Let $\inf_{x > K} {df \over dx}(x) > \bar{a} > 1$. It follows from $f(x) = f(K) + \int_{K}^x {df \over dx}(s) ds $ that for some $M < \infty$, $f(x) \geq M + \bar{a}x$ for $x > K$. Since both $q'_t$ and $m_t$ can take finitely many values, there exists $U$ such that $|u_t| \leq U$ for all $t$.  Let with $\gamma > 1$, a Lyapunov function be picked as $V(x)=\gamma^{-x}$, defined for positive $x$. Now, it follows that for sufficiently large $x$: $E[V(x_{t+1})|x_t=x] \leq V(x)$, since $E[\gamma^{-(f(x)+u_t+w_t)}] = E[\gamma^{-f(x)} \gamma^{-u_t} \gamma^{-w_t}] \leq E[\gamma^{-(M-U)} \gamma^{-(\bar{a} x+w_t})] = \gamma^{-(M-U)} \gamma^{-\bar{a} x} E[\gamma^{-w}]$, for all $x \in \{x : \gamma^{(\bar{a}-1)x} > E[\gamma^{-w}]\gamma^{-M+U)}\}$. Due to the additive noise process the source can escape any bounded interval with a non-zero probability. As a result, by Theorem 6.2.8 in \cite{YukselBasarBook} (see also Theorem 8.4.1 in \cite{MeynBook}), transience follows.\qed

Transience prohibits the existence of a stationary probability measure. The discussion above is parallel to Theorem 7.3.1 in \cite{YukselBasarBook} for linear systems. Related to the discussion above, for linear systems, the unboundedness of second moments in Proposition 5.1 in \cite{NairEvans} and the transience of such a controlled state process was established in Theorem 4.2 in \cite{YukBasTACNoisy}. We also note that \cite{Debasish} studied conditions for stabilization when the control actions are uniformly bounded, the controlled multi-dimensional system is marginally stable and is driven by noise with unbounded support. 

\subsection{Stationarity and Ergodicity}
In this section, instead of asymptotic mean stationarity, we will consider the more stringent condition of (asymptotic) stationarity of the controlled source process. For ease in presentation we will assume that $m_t$ takes values in a countable set, even though the extension to more general spaces is possible.

\begin{lem}
If the channel is memoryless, the process $(x_t, m_t)$ is a Markov chain.
\end{lem}

\textbf{Proof.} For any $t \in \mathbb{N}$,
\begin{eqnarray}
&& P(dx_t,m_t | x_s, m_s, s \leq t-1) \nonumber \\
&& = \sum P(dx_t,m_t, q'_{t-1} | x_s, m_s, s \leq t-1) \nonumber \\
&&= \sum P(dx_t| x_{t-1}, \gamma^d(m_{t-1},q'_{t-1})) P(q'_{t-1} |\gamma^e(x_{t-1},m_{t-1})) P(m_t | q'_{t-1},m_{t-1}) \nonumber \\
&& = \sum P(dx_t,m_t, q'_{t-1} | x_{t-1}, m_{t-1}) = P(dx_t,m_t | x_{t-1}, m_{t-1})
\end{eqnarray}
where we use the fact that the channel is of class A and (\ref{classAMarkov}) and (\ref{updates}). \qed

In the following, we assume that the channel is memoryless. For the Markov chain $(x_t, m_t)$, let $\pi_t(B) = P(x_t \in B)$ for all Borel $B$, that is, $\pi_t$ is the marginal occupation probability for the state process $x_t$. %Let $m_0$ be fixed in the following.

\begin{thm}\label{PHRCondition}
Suppose that the encoding, control and the memory update laws are given by (\ref{updates}). (i) Let $f$ have the form (\ref{model2}), (ii) Assumptions \ref{invertiblef} and \ref{bounds2} hold, (iii) $h(x_0) < \infty$. For the positive Harris recurrence of the process $x_t, m_t$ (which implies the existence of a unique invariant measure $\pi$ (and thus ergodicity)), it must be that
\begin{eqnarray}
C \geq  \int \pi(dx) \bigg(\int \nu(dw) \log_2(|J(f(x,w))|)\bigg), \label{boundI}
\end{eqnarray}
%provided that $(\int \nu(dw) \log_2(|J(f(x,w))|)$ is $\pi$-integrable and 
provided that $\limsup_{t \to \infty} {1 \over t} h(x_t) \leq 0$.
\end{thm}

\textbf{Proof.}
First note that
\begin{eqnarray*}
&& C \geq I(q_t,q'_t) = H(q'_t) - H(q'_t|q_t) \\
&& \geq H(q'_t | m_t) - H(q'_t|q_t) =  H(q'_t | m_t) - H(q'_t|q_t, x_t, m_t) \\
&& \geq H(q'_t | m_t) - H(q'_t|x_t, m_t) = I(q'_t;x_t|m_t)
\end{eqnarray*}
Hence,
\begin{eqnarray}
&&C \geq \liminf_{T \to \infty}  {1 \over T} \sum_{t=0}^{T-1} I(q'_t;x_t|m_t) \nonumber \\
&&= \liminf_{T \to \infty}  {1 \over T} \sum_{t=0}^{T-1} \bigg( h(x_t | m_t) - h(x_t | m_t,q'_t) \bigg) \nonumber \\
&&= \liminf_{T \to \infty}  \bigg( {1 \over T} \sum_{t=1}^{T-1} \bigg( h \bigg( f(x_{t-1},w_{t-1}) + Bu_{t-1} | m_t \bigg) - h(x_t | m_t,q'_t) \bigg) + I(q'_0;x_0|m_0) \bigg) \nonumber \\
&&= \liminf_{T \to \infty}  {1 \over T} \sum_{t=1}^{T-1} h\bigg( f(x_{t-1},w_{t-1}) + Bu_{t-1}) | m_t\bigg) - h(x_t | m_t,q'_t)  \nonumber \\
&&\geq  \liminf_{T \to \infty}  {1 \over T} \sum_{t=1}^{T-1} h\bigg(f(x_{t-1},w_{t-1}) + Bu_{t-1}| m_t, m_{t-1},q'_{t-1}\bigg) - h(x_t | m_t,q'_t)  \nonumber \\
&&=   \liminf_{T \to \infty}  {1 \over T} \sum_{t=1}^{T-1} h\bigg(f(x_{t-1},w_{t-1}) + Bu_{t-1}| m_{t-1},q'_{t-1}\bigg) - h(x_t | m_t,q'_t)  \nonumber \\
&&\geq \liminf_{T \to \infty}  {1 \over T} \sum_{t=1}^{T-1}  h\bigg( f(x_{t-1},w_{t-1}) + Bu_{t-1} | w_{t-1},m_{t-1},q'_{t-1}\bigg) - h(x_t | m_t,q'_t)  \nonumber \\
&&=  \liminf_{T \to \infty}  {1 \over T} \sum_{t=1}^{T-1} h\bigg(f(x_{t-1},w_{t-1}) | w_{t-1},m_{t-1},q'_{t-1}\bigg) - h(x_t | m_t,q'_t)  \nonumber \\
&&=  \liminf_{T \to \infty}  {1 \over T} \sum_{t=1}^{T-1} \bigg(\int \nu(dw) \bigg( \sum P(m_{t-1}=m, q'_{t-1}=q') \nonumber \\
&& \quad \quad \quad \quad \quad \quad \times \int P(x_{t-1} \in dx | m_{t-1}=m, q'_{t-1}=q',w_{t-1}=w) \log_2(|J(f(x,w_{t-1}))|) \nonumber \\
&& \quad \quad \quad \quad \quad \quad + h(x_{t-1}|m_{t-1}=m,q'_{t-1}=q',w_{t-1}=w) \bigg) - h(x_t | m_t,q'_t)\bigg) \nonumber \\
&&=  \liminf_{T \to \infty}  {1 \over T} \sum_{t=1}^{T-1} \bigg(\int \nu(dw) \bigg( \sum P(m_{t-1}=m, q'_{t-1}=q') \nonumber \\
&& \quad \quad \quad \quad \quad \quad \times \int P(x_{t-1} \in dx | m_{t-1}=m, q'_{t-1}=q') \log_2(|J(f(x,w_{t-1}))|) \nonumber \\
&& \quad \quad \quad \quad \quad \quad + h(x_{t-1}|m_{t-1}=m,q'_{t-1}=q') \bigg) - h(x_t | m_t,q'_t)\bigg) \label{indep} \\
&&=  \liminf_{T \to \infty}  {1 \over T} \sum_{t=1}^{T-1} \bigg(\int \nu(dw_{t-1}) \bigg(\int \pi_{t-1}(dx) \log_2(|J(f(x,w_{t-1}))|) \nonumber \\
&&  \quad \quad  \quad \quad  \quad \quad   \quad \quad  + h(x_{t-1}|m_{t-1},q'_{t-1}) \bigg)  - h(x_t  | m_t,q'_t) \bigg)  \label{equality}  \\
%&&=  \liminf_{T \to \infty}  {1 \over T}   \sum_{t=1}^{T-1} \int \nu(dw_{t-1}) \bigg(\int \pi_{t-1}(dx) \log_2(|J(f(x,w_{t-1}))|)  \nonumber \\
%&&  \quad \quad \quad  \quad \quad \quad  \quad \quad \quad  + h(x_{t-1}|m_{t-1},q'_{t-1}) - h(x_t  | m_t,q'_t) \bigg)   \nonumber \\
&&= \liminf_{T \to \infty}  {1 \over T}  \bigg(\sum_{t=1}^{T-1}  \int \pi_{t-1}(dx) \bigg(\int \nu(dw) \log_2(|J(f(x,w))|) \bigg) \nonumber \\
&&  \quad \quad  \quad \quad \quad  \quad \quad  \quad \quad   \quad \quad - h(x_{T-1} | m_{T-1},q'_{T-1})\bigg) \label{Fubinis} \\
&&\geq \liminf_{T \to \infty}  {1 \over T} \bigg( \sum_{t=1}^{T-1}  \int \pi_{t-1}(dx) \bigg(\int \nu(dw) \log_2(|J(f(x,w))|) \bigg) - h(x_{T-1}) \bigg) \nonumber \\
&&\geq \liminf_{T \to \infty}  {1 \over T}\sum_{t=1}^{T-1}  \int \pi_{t-1}(dx) \bigg(\int \nu(dw) \log_2(|J(f(x,w))|) \bigg) -  \limsup_{T \to \infty}  {1 \over T} h(x_{T-1}) \nonumber \\
%&&= \limsup_{t \to \infty} \int \pi_{t-1}(dx) (\int P(dw_t) \log_2(|J(f(x,w))|) )  + H(x_{t-1}|m_t)  - H(x_t | m_{t+1})  \\
&&\geq \liminf_{T \to \infty}  {1 \over T} \sum_{t=1}^{T-1} \int \pi_{t-1}(dx) \bigg( \int \nu(dw) \log_2(|J(f(x,w))|) \bigg) \label{SonDenklem}\\
&&= \liminf_{T \to \infty} \int \pi_0(dx) E_x\bigg[ {1 \over T} \sum_{t=0}^{T-1} \int \pi_{t-1}(dx) \bigg( \int \nu(dw) \log_2(|J(f(x,w))|) \bigg) \bigg] \label{Markov}\\
&&\geq \int \pi_0(dx)  \liminf_{T \to \infty} E_x\bigg[ {1 \over T} \sum_{t=0}^{T-1} \int \pi_{t-1}(dx) \bigg( \int \nu(dw) \log_2(|J(f(x,w))|) \bigg) \bigg] \label{Fatou}\\
&&= \int \pi_0(dx) \bigg( \int \pi(dz) \bigg( \int \nu(dw) \log_2(|J(f(z,w))|) \bigg) \bigg) \label{fnorm}\\
&& = \int \pi(dx) \bigg(\int \nu(dw) \log_2(|J(f(x,w))|) \bigg)
 \label{SonDenklem}
\end{eqnarray}
In the first lines above, we use the fact that conditioning on a random variable reduces the entropy and the update laws (\ref{updates}). The equality (\ref{equality}) holds since for every $w$, the map $f(.,w)$ is invertible (here $J(f(x,w_t))$ is the Jacobian for the realized value of $w_t$) and that $w_t$ is an independent noise process using the laws of total probability. Here (\ref{indep}) follows due to the independence of $w_t$, (\ref{Fubinis}) follows from Fubini's Theorem since $\log_2(|J(f(x,w))|)$ is bounded, (\ref{Fatou}) follows from Fatou's lemma given the assumption that $ \log_2(|J(f(x,w))|)$ is bounded from below, and (\ref{fnorm}) follows from positive Harris recurrence (see \cite[Theorem 4.3.1]{HernandezLasserreErgodic}).
\qed

\begin{remark}
If one considers a more general control-affine model such as of the form (\ref{model1}) with $x_{t+1}=f(x_t,w_t)+B(u_t)x_t$, the condition would read as:
\[C \geq \int \pi(dx,m,q') \bigg(\int \nu(dw) \log_2\bigg( \bigg|J\bigg(f(x,w) + B(\gamma^d(m,q'))x \bigg) \bigg| \bigg) \bigg),\]
where $\pi$ is invariant for the (enlarged) Markov chain $(x_t,m_t,q'_t)$. \hfill $\diamond$ 
\end{remark}

\section{Discrete Noiseless Channels and a Stationary and Ergodic Construction}\label{Sec3}

In this section, we provide achievability results and a stabilizing coding/control policy.  As discussed earlier, the study of non-linear systems have typically considered noise-free controlled systems; e.g. \cite{baillieul2004data}, \cite{liberzon2005stabilization}, and \cite{de2004stabilizability}. As also noted earlier, for noise-free systems, it typically suffices to only consider a sufficiently small {\it invariant} neighborhood of an equilibrium point to obtain stabilizability conditions which is not necessarily the case when the system is driven by an additive noise process. We consider such an example in the following. 

\begin{thm}\label{theoremAchieve2}
Consider a non-linear system of the form (\ref{model5}), where $\{w_t\}$ is a sequence of zero-mean Gaussian random vectors and there exists a control function $\kappa(z)$ such that $|f(x, \kappa(z))|_{\infty} \leq |a| |x - z|_{\infty}$ for all $x, z \in \mathbb{R}^N$, with $\kappa(0)=0$. 
% and $|f(x,0)|_{\infty} \leq |a| |x|_{\infty}$. 
%Let $f$ be Lipschitz continuous with the $\sup$-norm coefficent:
%\[\sup_{x,y} {|f(x) - f(y)|_{\infty} \over |x-y|_{\infty}} = a > 1\]
For the stationarity and ergodicity of $\{x_t\}$ (and thus with a unique invariant probability measure), it suffices that $C >N\log_2(|a|) + 1$. 
\end{thm}

\begin{remark}
It may be possible in general to reduce the rate requirements by the use of variable-rate encoding schemes; for example, if there exists a compact region outside of which the constant $a$ can be upper bounded by a smaller number, a region-dependent quantization rate can be applied which can reduce the average data rate required for system stability. In this paper, since there is an explicit channel, our focus has been on fixed-rate coding schemes. \hfill $\diamond$
\end{remark}

%\begin{remark}
%Consider a special case of the non-linear system of the form (\ref{model3}), where $x_t$ is an $\mathbb{R}^n$-valued random variable: $x_{t+1}=f(x_t)+Bu_t + w_t$. If $f$ be Lipschitz continuous with the $\sup$-norm coefficent:
%\[\sup_{x,y} {|f(x) - f(y)|_{\infty} \over |x-y|_{\infty}} = a > 1,\] $B$ is full-rank and $f(0)=0$, Theorem \ref{theoremAchieve2} holds. 
%%\begin{thm}\label{theoremAchieve}
%%Let $f$ be Lipschitz continuous with the $\sup$-norm coefficent:
%%\[\sup_{x,y} {|f(x) - f(y)|_{\infty} \over |x-y|_{\infty}} = a > 1\]
%%For the existence of a unique invariant probability measure, it suffices that $C > (n \log_2(|a|) + 1)$. For the AMS property, it suffices that $C > n \log_2(|a|)$.
%%\end{thm}
%%
%%\begin{remark}
%%We note that the achievability conditions on the system can be relaxed to a setup of the form (\ref{model4}), with the assumption of the existence of a control function $\kappa(z)$ such that $|f(x, \kappa(z)),w|_{\infty} \leq |a| |x - \hat{x}|_{\infty} + c|w|_{\infty}$, for some $c < \infty$, $\kappa(0)=0$ and $|f(x,0)|_{\infty} \leq |a| |x|_{\infty}$. 
%%\end{remark}
%\end{remark}

\textbf{Proof.} 
The proof follows essentially from the approach developed in \cite{YukTAC2010} and \cite{AndrewJohnstonReport} with extension to non-linear analysis. Consider the case with $N=2$. Let $\Delta >0$ denote the bin size for a uniform quantizer and let for each coordinate $x^i \in \mathbb{R}, i=1,2;$
\begin{eqnarray}
Q_K^{\Delta}(x^i) = \begin{cases}   (k - \half (K + 1) ) \Delta,  \quad \quad &
	\mbox{if} \ \ x^i \in [ (k-1-\half K   ) \Delta , (k-\half K  ) \Delta)  \\
 (\half (K - 1) ) \Delta,  \quad \quad &
	\mbox{if} \ \ x^i = \half K \Delta \\
  0 ,  \quad \quad &
  	\mbox{if} \ \ x^i \not\in [- \half K  \Delta, \half K  \Delta].
\end{cases} \label{UniformQuantizerDescriptionCHP8}
\end{eqnarray}
and define ${\bf Q}^{\Delta}(x)= (Q_K^{\Delta}(x^1),Q_K^{\Delta}(x^2))$ if $Q_K^{\Delta}(x^i) \neq 0$ for $i=1,2$ and ${\bf Q}^{\Delta}(x)=0$ if $Q_K^{\Delta}(x^i) = 0$ for some $i$. Thus, the number of symbols in the image of ${\bf Q}^{\Delta}$ is $K^2+1$ (and not $(K+1)^2$). The quantizer outputs are transmitted through a memoryless erasure channel, after being subjected to a bijective mapping, which is performed by the channel encoder. The channel encoder maps the quantizer output symbols to corresponding channel inputs $q \in \clM\eqdef\{1,2\dots,K^2+1\}$. A channel encoder at time $t$, denoted here by ${\cal E}_t$, maps the quantizer outputs to $\clM$ such that ${\cal E}_t(Q_t(x_t))=q_t \in \clM$. For $i=1,2$, let $R'=\log_2(K)$. For $t \geq 0$ and with $\Delta^1_0 = \Delta^2_0 \in \mathbb{R}$, define
\[h^i_t = { x^i_t \over \Delta^i_t 2^{R'-1}},\]
and with
\[\hat{x}_{t} = \begin{bmatrix}\hat{x}^1_t \\ \hat{x}^2_t \end{bmatrix},\]
consider:
\begin{eqnarray}
\label{QuantizerUpdate4CHP7}
&&u_t = - \kappa(\hat{x}_{t}), \nonumber \\
&& \begin{bmatrix}\hat{x}^1_t \\ \hat{x}^2_t \end{bmatrix} = \begin{bmatrix} Q_{K_1}^{\Delta^1_t}(x^1_t)  \\ Q_{K_2}^{\Delta^2_t}(x^2_t) \end{bmatrix} 1_{\{\max_i |h_i| \leq 1\}} + \begin{bmatrix} 0  \\ 0 \end{bmatrix} 1_{\{\max_i |h_i| > 1\}}, \\
&&\Delta^1_{t+1} = \Delta^1_t \bar{Q}(|h^1_t|,|h^2_t|,\Delta^1_t, \Delta^2_t), \quad \quad \Delta^2_{t+1} = \Delta^2_t \bar{Q}(|h^1_t|,|h^2_t|,\Delta^1_t, \Delta^2_t),
\end{eqnarray}
with, for $i=1,2$, $\delta > 0$ $\alpha \in (0,1)$, $L > 0$ such that
\begin{eqnarray}
 \bar{Q}(x,y,\Delta^1, \Delta^2) &=& |a| + \delta \quad \mbox{if } \quad |x| > 1, \quad \mbox{or} \quad \quad |y| > 1  \nonumber \\
  \bar{Q}(x,y,\Delta^1, \Delta^2) &=& \alpha \quad \mbox{if } \quad |x| \leq 1, |y| \leq 1; \quad \Delta^1 > L, \Delta^2 > L   \nonumber \\
 \bar{Q}(x,y,\Delta^1, \Delta^2) &=& 1 \quad \quad  \mbox{if } \quad |x| \leq 1, |y| \leq 1; \quad \Delta^1 \leq L \quad \mbox{or} \quad  \Delta^2 \leq L  \nonumber
\end{eqnarray}
%with $\sqrt {\Expect[w_t^2]} < \delta L{|a| \over 2^{R'} - \eta } 2^{R'-1}.$
Note that, the above imply $\Delta^i_t \geq \alpha L$. To make the state space for the bin size process countable as in \cite{YukselBasarBook} \cite{YukTAC2010}, we take that $\log_2(\bar{Q}(\cdot))$ take values in integer multiples of $s$ where the integers taken are relatively prime (that is they share no common divisors except for $1$); see \cite[Lemma 7.6.2]{YukselBasarBook}. 

%\begin{figure}[h]
%\centering
%\epsfig{figure=ModifiedUniformVectorQuantizer.eps,height=7cm,width=11cm}
%%\includegraphic{figure=truncpdf.ps}
%\caption{The vector quantizer. There is a single overflow bin. \label{VectorFigureQuantizer}}
%\end{figure}

We note the following without proof. 
\begin{lem} \label{MarkovChain}
The process $(x_t,\Delta_t)$ is a Markov chain.
\end{lem}

We define a sequence of stopping times as follows:
\begin{eqnarray} \label{vectorStoppingTimes}
\stp_0 = 0, \quad \stp_{z+1} &=& \inf \{k > \stp_z : |h^i_{k}| \leq  1, i \in \{1,2\} \}, \quad z \in \mathbb{Z}_+. \nonumber
\end{eqnarray}

By the strong Markov property and the nature of the stopping times, $(x_{\stp_z}, h_{\stp_z})$ is also Markov. In the following, we show that there exist $b_0 >0$, $b_1 < \infty$ such that
\begin{eqnarray}\label{geoDrift}
\Expect[\log(\Delta_{\stp_{z+1}}^2)| \Delta_{\stp_z }, h_{\stp_z }] \leq \log(\Delta_{\stp_z }^2) - b_0 + b_1 1_{\{|\Delta_{\stp_z}| \leq F\}}
\end{eqnarray}

We first bound the probability $\Prob(\stp_{z+1} - \stp_z \geq k|\Delta_{\stp_z },h_{\stp_z })$ from above.
\begin{lem}
 The discrete probability measure $\Prob(\stp_{z+1}-\stp_z=k\mid x_{\stp_z },\Delta_{\stp_z })$ has the upper bound
\[
\Prob(\stp_{z+1}-\stp_z \geq k | x_{\stp_z },\Delta_{\stp_z }) \leq M(\Delta_{\stp_z}) r^{-k},
\]
for some $r > 1$ and $\lim_{\Delta \to \infty}M(\Delta) = 0$. 
\end{lem}

\textbf{Proof.}
Observe that for $0 < k < \tau_1$, $\hat{x}_k= 0$ and $u_k= \kappa(0)=0$. Let $|x| = \|x\|_{\infty}$. Now, for $k \geq 2$
\begin{eqnarray}
&& \Prob(\stp_{1} \geq k| x_0, \Delta_0) \leq \Prob_{x_0,\Delta_{0}}\bigg( |x_{k-1}| \geq (|a|+\delta)^{k-2}2^{R'-1} \alpha \Delta_0  \bigg) \nonumber \\
&& \le \Prob_{x_0,\Delta_{0}}\bigg( |f(x_{k-2})| + |w_{k-2}| \geq (|a|+\delta)^{k-2}2^{R'-1} \alpha \Delta_0  \bigg) \nonumber \\
&& \le \Prob_{x_0,\Delta_{0}}\bigg( |a(x_{k-2})| + |w_{k-2}| \geq (|a|+\delta)^{k-2}2^{R'-1} \alpha \Delta_0  \bigg) \nonumber \\
&& \le \Prob_{x_0,\Delta_{0}}\bigg( |a| |x_{k-2}| + |w_{k-2}| \geq (|a|+\delta)^{k-2}2^{R'-1} \alpha \Delta_0  \bigg) \nonumber \\
&& \le \Prob_{x_0,\Delta_{0}}\bigg( \sum_{i=0}^{k-2} |a|^{-i} |w_{i}| \geq {(|a|+\delta)^{k-2}2^{R'-1} \alpha \Delta_0 \over |a|^{k-1}} - |x_0 - \hat{x}_0| \bigg) \label{inductiveStep} \\
&& \le \Prob_{x_0,\Delta_{0}}\bigg( \sum_{i=0}^{k-2} |a|^{-i} |w_{i}| \geq {(|a|+\delta)^{k-2}2^{R'-1} \alpha \Delta_0 \over |a|^{k-1}} - \Delta_0/2 \bigg) \nonumber \\
&& = \Prob_{x_0,\Delta_{0}}\bigg( \sum_{i=0}^{k-2} |a|^{-i} |w_{i}|  \geq \Delta_0/2 \bigg( ({|a|+\delta \over |a|})^{k-2}{2^{R'} \alpha \over |a|}  - 1 \bigg) \bigg) \label{comp} \\
&& \leq { E[ \sum_{i=0}^{\infty} |a|^{-i} |w_{i}|] \over  \Delta_0/2 \bigg(({ |a|+\delta \over |a|})^{k-2}{2^{R'} \alpha \over |a|}  - 1 \bigg) } \label{comp1} \\
&& \leq M(\Delta_0) r^{-k}  \label{comp2}
\end{eqnarray}
with \[M(\Delta_0) = {K E[ \sum_{i=1}^{\infty} |a|^{-i} |w_{i}|] \over ({|a|+\delta \over |a|})^2 \Delta_0 ({2^{R'} \alpha \over |a|}  - 1)} < \infty,\]
for some $K < \infty$ and $r \in (1,(|a|+\delta)/|a|)$ so that $\lim_{\Delta_0 \to \infty} M(\Delta_0)= 0$. Here, (\ref{inductiveStep}) follows from an inductive argument, (\ref{comp}) follows from the fact that the term \[\bigg(2^{R'-1}({|a|+\delta \over |a|})^{k-2} {\alpha \over |a|} - {1 \over 2} \bigg)\]
 is positive for $k \geq 2$ provided that $2^{R'} > {|a| \over \alpha}$, (\ref{comp1}) follows from Markov's inequality and (\ref{comp2}) from the fact that $w_i$ is Gaussian together with the property $|w| \leq N(1 + |w|^2)$ leading to the finiteness of $E[ \sum_{i=1}^{\infty} |a|^{-i} |w_{i}|]$.  \hfill $\diamond$

%NOTE: CHECK THE COEFFICIENT AGAIN GIVEN THE NEW MORE GENERAL MODEL...CHECK...

We now invoke \cite[Theorem 2.1]{YukMeynTAC2010}: Let $\bfmX$ be an $\mathbb{X}$-valued Markov chain (where $\mathbb{X}$ is a standard Borel space) and $\stp_z, z \geq 0$ be a sequence of stopping times measurable on the filtration generated by the state process with $\stp_0=0$.
\begin{thm}
\label{thm5} \cite[Theorem 2.1]{YukMeynTAC2010}
Suppose that $\bfmX$ is a $\varphi$-irreducible and aperiodic Markov chain. Suppose moreover that there are functions
 $V \colon \mathbb{X} \to (0,\infty)$,
 $\delta \colon \mathbb{X} \to [1,\infty)$,
 $f\colon \mathbb{X} \to [1,\infty)$,
%\notes{Must make $f\ge 1$ to ensure positivity}
a small set $C$ on which $V$ is bounded, and a constant $b \in \Re$,   such that the following hold:
%\notes{Note changes to make similar to previous theorem.  And note I removed kappa -- better check proofs}
\begin{equation}
\begin{aligned}
\Expect[V(x_{\stp_{z+1}}) \mid \clF_{\stp_z }]  &\leq  V(x_{\stp_{z}}) -\delta(x_{\stp_z }) + b1_{\{x_{\stp_z} \in C\}}
\\
 \Expect \Bigl[\sum_{k=\stp_z}^{\stp_{z+1}-1} f(x_k)  \mid \clF_{\stp_z }\Bigr]  &\le \delta(x_{\stp_z})\, , \qquad \qquad \qquad \qquad z\ge 0.
\end{aligned}
\label{e:thm5delta}
\end{equation}
Then the following hold:
\begin{romannum}
\item
$\bfmX$ is positive Harris recurrent, with unique invariant distribution $\pi$
\item $\pi(f)\eqdef \int f(x)\, \pi(dx) <\infty$
\item  For any function $g$ that is bounded by $f$, in the sense that $\sup_{x} |g(x)|/f(x)<\infty$, we have convergence in the mean, and the Law of Large Numbers holds:
%\archive{To do:  Add proof.
%\\
%I added this result so the reader can see some motivation}
\[
\begin{aligned}
\lim_{t\to\infty} \Expect_{x}[g(x_t)] &= \pi(g)
\\
\lim_{N\to\infty} \frac{1}{N} \sum_{t=0}^{N-1} g(x_t) &= \pi(g)\qquad a.s.\,, \ x\in \mathbb{X}
\end{aligned}
\]
\end{romannum}
\end{thm}

By taking $f(x)=1$ for all $x \in \mathbb{X}$, the following holds.
\begin{thm}\label{corol}\cite{YukMeynTAC2010}
Suppose that $\bfmX$ is a $\varphi$-irreducible Markov chain with natural filtration ${\cal F}_t$.   Suppose moreover that there is a function $V: \mathbb{X} \to (0,\infty)$, a small set $C$ on which $V$ is bounded, and a constant $b \in \Re$,   such that the following hold:
%\notes{note:  crucial to sup over z}
%\[
\begin{equation}\label{PositiveRecur}
\begin{aligned}
\Expect[V(x_{\stp_{z+1}}) &\mid \clF_{\stp_z }] \leq V(x_{\stp_z }) - 1 + b1_{\{x_{\stp_z} \in C\}}
\\
\sup_{z\ge 0} \Expect[\stp_{z+1} -\stp_z &\mid \clF_{\stp_z }] < \infty.
\end{aligned}
\end{equation}
%\]
Then $\bfmX$ is positive Harris recurrent.
\end{thm}

Now, with the candidate Lyapunov function $V_0(x_t,\Delta_t) = \log(\Delta^2_t)$, for $\Delta_{\stp_z }> L$,
\begin{eqnarray}
&& \Expect[V_0(x_{\stp_{z+1}},\Delta_{\stp_{z+1}}) \mid x_{\stp_z},\Delta_{\stp_z}] = \Prob(\stp_{z+1} - \stp_z = 1) \bigg(2 \log(\alpha) + \log(\Delta_{\stp_z }^2) \bigg) \nonumber \\
&&\quad \quad \quad \quad \quad \quad \quad + \sum_{k=2}^{\infty} \log(\Delta_{\stp_{z}+k}^2) \Prob(\stp_{z+1} - \stp_z = k\mid x_{\stp_z},\Delta_{\stp_z}) \nonumber \\
&&=  \Prob(\stp_{z+1} - \stp_z = 1) \bigg(2 \log(\alpha) + \log(\Delta_{\stp_z }^2) \bigg) \nonumber \\
&&\quad \quad \quad \quad \quad \quad \quad + \sum_{k=2}^{\infty} 2(\log_2(\alpha)+ (k-1) (|a|+\delta) M(\Delta) r^{-k} \nonumber
%&& \quad\quad\quad + \log(\Delta_{\stp_z}^2) + 2\log({|a| \over 2^{R'} - \eta})  \Prob(\stp_{z+1} - \stp_z \geq 2|\Delta_{\stp_z },h_{\stp_z %})\nonumber \\
%&& = \bigg(2\log({|a| \over 2^{R'} - \eta}) + \log(\Delta_{\stp_z }^2) \bigg) + \sum_{k=2}^{\infty} \log(\Delta_{\stp_{z}+k}^2) \Prob(\stp_{z+1} - \stp_z = k\mid x_{\stp_z},\Delta_{\stp_z}) \nonumber
\end{eqnarray}

Now, by (\ref{comp2}), $\lim_{\Delta_0 \to \infty} \Prob(\stp_{z+1} - \stp_z = 1 | \Delta_0,x_0) = 1$ uniformly in $|x_0| \leq 2^{R'-1}\Delta_0$.
As a result, the drift condition of Theorem \ref{corol} holds. We need to ensure, however, the small/petite set \cite{MeynBook} property of compact sets to establish positive Harris recurrence. A sufficiently small compact set for this chain is petite due to the countability of the values that $\Delta_t$ takes and the {\it uniform countable additivity} property of the Markov chain due to the presence of the additive Gaussian noise, as in p.206 of \cite{YukselBasarBook} and the continuity of $f$ in $x$. This argument applies for $N$-dimensional systems as well with $N > 2$. This completes the proof of Theorem \ref{theoremAchieve2}.
\qed

% or a quasi-continuity argument NOTE: DO we need to state the continuity of f, we already assume continuity in $x,u$ so it should be sufficient. \qed

\begin{remark}
The approach adopted in the proof of Theorem \ref{theoremAchieve2} applies for more general channels (such as erasure channels or discrete memoryless channels) subject to more tedious error bounds. \hfill $\diamond$ %One could also apply the same approach to establish stronger versions of stability, such as quadratic stability.  \hfill $\diamond$
\end{remark}

\section{Discussion and conclusion}\label{Sec4}

In this paper, conditions on information channels leading to the stochastic stability of non-linear systems controlled over noisy channels has been investigated. Stochastic stability notions considered were asymptotic mean stationarity, ergodicity and stationarity. Results for linear systems are recovered as a special case.

In the following we present some future directions and a comparison with the results involving topological entropy.

\subsection{Comparison with invariance entropy and deterministic non-linear systems controlled over noiseless channels}

As noted earlier, noise-free systems and noiseless discrete channels have been studied in the literature in the context of topological entropy and invariance entropy. Here, we establish some connections.

One related result in this literature is with regard to stabilization to a point: Under the assumptions that (i) $f$ has the form in (\ref{model5}) (without noise) with continuous partial derivatives, (ii) there exists a fixed point (equilibrium) $x^*$ so that $x^*=f(x^*,u^*)$, (iii) a local strong invariability condition is satisfied which relates the size of an invariant set and the size of a control action set in the sense that for any $ \epsilon >0$, there exist $\rho > 0$ so that for all $\epsilon' \in (0,\rho]$, the set $\{x: |x-x^*| \leq \epsilon'\}$ is strongly invariant with the control action set $U = \{u: |u-u^*| \leq \epsilon\}$, and (iv) the pair $(A,B)$ is controllable where $A, B$ are the Jacobians of $f$ with respect to state and control at $x^*,u^*$, \cite{nair2004topological} has reported that for convergence to the equilibrium an average rate $R > \sum_{|\lambda_i| > 1} \log_2(|\lambda_i|)$
 is sufficient, where $\lambda_i$ are the eigenvalues of the Jacobian at the equilibrium point. 
 
A further related result in spirit to our paper is on a case where there exists an invariant set with a non-empy interior: For continuous-time systems of the form ${dx \over dt} = f(x,u), u \in U$, Colonius and Kawan \cite{colonius2009invariance} establish a lower bound on {\it invariance entropy} as
\begin{eqnarray}\label{KCBound}
\max\bigg(0,\min_{x,u \in Q \times U} \sum_i {\partial f_i \over \partial x_i}(x,u) \bigg),
\end{eqnarray}
 where $Q$ is a weakly invariant set and $f_i$ is the $i$th coordinate function of $f$. More refined bounds are present if further structural properties are imposed: in \cite{KawanDCDS2016}, under a uniform hyperbolicity assumption (see \cite[Definition 4.4]{KawanDCDS2016}), Theorem 4.8 states a similar lower bound by considering the unstable components in an invariant set. 

These results can be viewed to be related to Theorem \ref{notAMS} and Theorem \ref{PHRCondition}, as well as Theorem \ref{theoremAchieve2}, in that the average entropy growth as measured by the eigenvalues of the Jacobian matrix under the invariant probability measure is lower bounded by a minimum over the elements in the support set, or is upper bounded by a maximizing element in the support set. In the stabilization to the point example of \cite{nair2004topological}, the invariant measure is a delta measure on a single point. In the invariant set example leading to (\ref{KCBound}), the set $Q$ can be viewed to be the support set of some invariant measure under the system dynamics if such a measure were to exist. 
Likewise, \cite{liberzon2005stabilization} and \cite{liberzon2002stabilization} have obtained conditions for noise-free systems controlled over noiseless channels. Due to the absence of noise, one could identify an invariant compact set, and consider a bound on the Lipschitz growth parameter for the system over this invariant set to obtain sufficiency conditions. When the system is (Lebesgue) irreducible, however, due to the effect of noise, local properties are not descriptive and the invariant probability measure reflects the rate conditions and entropy growth in the system. In this case, the local growth integrated under an invariant measure gives a proper bound.

Differential entropy is a useful measure for how much a stochastic system generates uncertainty, however our analysis does not distinguish between the stable and unstable modes of a controlled system and is only able to resemble the classical results in ergodic theory (Pesin's formula \cite{young2003entropy}) for {\it expanding} systems, and thus, with only positive Lyapunov exponents. In the linear case, the arguments follow by restricting the state space to those corresponding to the unstable modes. For a general non-linear system, however, a careful geometric study needs to be done. On the other hand, for deterministic systems, under a topological entropy formulation, the rate of growth can be measured by local Jacobian matrices, but such a topological discussion requires further geometric analysis with regard to the use of appropriate metrics, as studied extensively in \cite{kawan2013invariance}. Thus, the connection between the differential entropy method and geometric approaches requires some further study.

We note also that recently a metric entropy generalization of some of the results in \cite{kawan2013invariance} have been developed \cite{Colonius2016}. 

\subsection{Some open directions on stationary coding and control policies and information theory}
It would be interesting to show, for a class of systems, that stationary coding and control policies can be used to arrive at stability with a stationary closed loop-process provided that the capacity of the channel satisfies the entropy growth bound and the channel satisfies certain ergodicity conditions. However, except for linear Gaussian systems controlled over Gaussian channels and erasure channels (see \cite{YukselBasarBook} for a detailed discussion for both setups), this question has not been answered even for linear systems controlled over general discrete memoryless channels (that is, non-stationary coding schemes have been used for more general discrete memoryless channels). Furthermore, the tightness of the converse results is another direction. A further direction is the causal coding problem for non-ergodic sources: In the information theory literature, through non-causal codes, a class of source coding theorems for non-ergodic sources exist (see e.g. \cite{gray1974source}), however, the extensions of these for even control-free non-linear systems under causal coding require further research.

\section{Acknowledgements}

We gratefully acknowledge extensive technical discussions with Dr. Christoph Kawan and the suggestions of an anonymous reviewer.

\appendix

\section{Stationary, ergodic, and asymptotically mean stationary processes} \label{ergodicTheory}
\index{Ergodic processes}
In this subsection, we review ergodic theory, in the context of information theory (that is with the transformations being specific to the shift operation). A comprehensive discussion is available in Shields \cite{Shields}, Gray \cite{GrayProbabilit}, \cite{GrayKieffer}, and Appendix C in \cite{YukselBasarBook}. 

Let $\mathbb{X}$ be a complete, separable, metric space. Let ${\cal B}(\mathbb{X})$ denote the Borel sigma-field of subsets of $\mathbb{X}$. Let $\Sigma=\mathbb{X}^{\infty}$ denote the sequence space of all one-sided or two-sided infinite sequences drawn from $\mathbb{X}$. Thus, for a two-sided sequence space if $x \in \Sigma$ then $x = \{\dots,x_{-1},x_0,x_1,\dots\}$ with $x_i \in \mathbb{X}$. Let $X_n:\Sigma \to \mathbb{X}$ denote the coordinate function such that $X_n(x)=x_n$. Let $T$ denote the shift operation on $\Sigma$, that is $X_n(Tx)=x_{n+1}$. That is, for a one-sided sequence space $T(x_0,x_1,x_2,\dots)=(x_1,x_2,x_3,\dots)$.

Let ${\cal B}(\Sigma)$ denote the smallest sigma-field containing all cylinder sets of the form $\{x: x_i \in B_i, m \leq i \leq n\}$ where $B_i \in {\cal B}(\mathbb{X})$, for all integers $m,n$. Observe that $\cap_{n \geq 0} T^{-n}{{\cal B}(\Sigma)}$ is the tail $\sigma$-field $\cap_{n \geq 0} \sigma(x_n,x_{n+1},\cdots)$, since $T^{-n}(A) = \{x: T^nx \in A\}$.

Let $\mu$ be a {\it stationary} measure on $(\Sigma,{\cal B}(\Sigma))$ in the sense that $\mu(T^{-1}B)=\mu(B)$ for all $B \in {\cal B}(\Sigma)$. Then, the sequence of random variables $\{x_n\}$ defined on the probability space $(\Sigma,{\cal B}(\Sigma),\mu)$ is a stationary process.

\begin{defn} Let $\mu$ be the measure on a process. This random process is ergodic if $A=T^{-1}A$ implies that $\mu(A) \in \{0,1\}$.
\end{defn}

%A source is {\em ergodic} if $TB=B$ implies that $\mu(B) \in \{0,1\}$,
That is, the events that are unchanged with a shift operation are trivial events.
%All other events, in an ergodic process, are aperiodic.
%If a process is ergodic then the following holds (IS THE FOLLOWING TRUE FOR ALL $x$ VALUES, OR ONES WHICH HAS A LIMIT AND IN AN INVARIANT SET?, BUT ALL SETS ARE INVARIANT, SO OK): Let $f$ be %integrable. Then, there exists an integrable function $\hat{f}$ such that:

Mixing is a sufficient condition for ergodicity. Thus, a source is ergodic if $\lim_{n \to \infty} P(A \cap T^{-n}B) = P(A) P(B)$, since the process forgets its initial condition. For the special case of Markov sources, we have the following: A positive Harris recurrent Markov chain is ergodic, since such a process is mixing and stationary.
%A random process is $N-$ergodic if it is ergodic with respect to $T^N$ \cite{GrayProbabilit}. We will be primarily concerned with such processes in %this paper.
%
%\begin{defn} A random process with measure $\mu$ is $N$-stationary, (cyclo-stationary or periodically stationary with period $N$) if $\mu(T^{-N}B)=\mu(B)$ for all $B \in {\cal B}(\Sigma)$, or equivalently for any $n \in \mathbb{N}$ samples $t_1, t_2, \dots, t_n$:
%\[\mu(x_{t_1} \in A_1, \dots,x_{t_n} \in A_n) = \mu(x_{t_1 + N} \in A_1, \dots,x_{t_n + N} \in A_n) \]
%\end{defn}
%\begin{defn} A random process is $N$-ergodic if $A=T^{-N}A$ implies that $\mu(A) \in \{0,1\}$.
%\end{defn}
%
%%\begin{defn}\label{CoordinateRecurrent} A set $A \in {\cal B}(\mathbb{X})$ is coordinate-recurrent if for some $m \in \mathbb{Z}_+$
%%\[\sum_{m=0}^{\infty} 1_{\{X_m(x) \in A \}} = \infty, \quad a.s.\]
%%\end{defn}
%%%

\begin{defn}\index{Asymptotic mean stationarity (AMS)} \label{AMSDefinition}
A process on a probability space $(\Omega, {\cal F}, {\bf P})$ with process measure $P$, is asymptotically mean stationary (AMS) if there exists a probability measure $\bar{P}$ such that
\[\lim_{N \to \infty} {1 \over N} \sum_{k=0}^{N-1} P(T^{-k} F) = \bar{P} (F),\]
for all events $F \in {\cal B}(\Sigma)$. Here $\bar{P}$ is called the stationary mean of $P$, and is a stationary measure.
\end{defn}
Note that ${\bar P}$ is stationary since, by definition $\bar{P}(F) = \bar{P}(T^{-1}F)$. For the importance of the AMS property, its relations with Birkhoff's ergodic theorem, some applications and sufficient conditions, please see \cite{GrayProbabilit} or \cite{GrayKieffer}.%A cyclo-stationary process is AMS. 

\section{Proof of Theorem \ref{ModifiedMatveevsLemma}}\label{ProofofModified}

%\textbf{Proof of Proposition \ref{ModifiedMatveevsLemma}.}
Define the event for $K > 0$ so that $P(|x_0| < K) > 0$ as
\[{\cal S}^K_{\eta} = \{\omega: |x_0| \leq K, w=\eta, i.e., w_k = \eta_k, \eta_k \in \mathbb{R}^p, k \geq 0\},\] 
such that the noise realizations are fixed and deterministic. In the following, we will drop the subscript and superscripts and let $P_{{\cal S}}$ or $P(\cdot | {\cal S})$ denote the conditional probabilities given the event ${\cal S}^K_{\eta}$.We recall here that $\{w_t, t \geq 0\}$ and $x_0$ are assumed to be independent. By Definition \ref{ClassAChannelDefinition}, first note that the capacity expression satisfies
\begin{eqnarray}
C&=& \lim_{T \to \infty} \max_{\{ P(q_t | q_{[0,t-1]},q'_{[0,t-1]}), \quad 0 \leq t \leq T-1\}} {1 \over T} I(q_{[0,T-1]} \to q'_{[0,T-1]}) \nonumber \\
&=& \lim_{T \to \infty} \max_{\{ P(q_t | q_{[0,t-1]},q'_{[0,t-1]}), \quad 0 \leq t \leq T-1\}} {1 \over T} I(q_{[0,T-1]} \to q'_{[0,T-1]} | {\cal S}), \label{EqCon}
\end{eqnarray}
where the conditional directed information is given by \[I(q_{[0,T-1]} \to q'_{[0,T-1]} | {\cal S}) = \sum_{t=1}^{T-1} I(q_{[0,t]};q'_t|q'_{[0,t-1]}, {\cal S}) + I(q_0;q'_0|{\cal S}).\]
Here, (\ref{EqCon}) is a result of the following: Consider an encoder policy given by \[P^*=\{P^*(q_0), P^*(q_1|q_0,q'_0), \cdots, P^*(q_t|q_{[0,t-1]},q'_{[0,t-1]}), \cdots\}.\] For any $t \in \mathbb{N}$, almost surely the following holds:
\begin{eqnarray}
&& P(q'_t| q'_{[0,t-1]}, {\cal S})\nonumber \\
&& = \sum_{q_{[0,t]}} P(q'_t, q_{[0,t]}| q'_{[0,t-1]}, {\cal S})\nonumber \\
&& = \sum_{q_{[0,t]}} P(q'_t | q_{[0,t]}, q'_{[0,t-1]}, {\cal S}) P(q_{[0,t]}| q'_{[0,t-1]}, {\cal S})\nonumber \\
&& = \sum_{q_{[0,t]}} P(q'_t | q_{[0,t]}, q'_{[0,t-1]}) P(q_{[0,t]}| q'_{[0,t-1]}, {\cal S}) \label{defn00} \\
&& = \sum_{q_{[0,t]}} P(q'_t | q_{[0,t]}, q'_{[0,t-1]}) P^*(q_t| q'_{[0,t-1]}, q_{[0,t-1]}, {\cal S}) P(q_{[0,t-1]} | q'_{[0,t-1]}, {\cal S})\nonumber \\
&& = \sum_{q_{[0,t]}} P(q'_t | q_{[0,t]}, q'_{[0,t-1]}) P^*(q_t| q'_{[0,t-1]}, q_{[0,t-1]}) P(q_{[0,t-1]} | q'_{[0,t-1]}, {\cal S})\label{capD} \\
 && = \sum_{q_{[0,t]}} P(q'_t | q_{[0,t]}, q'_{[0,t-1]}) P^*(q_t| q'_{[0,t-1]}, q_{[0,t-1]}) P(q_{[0,t-1]} | q'_{[0,t-1]}) \label{sonD}\\
% {P(q_{[0,t-1]},q'_{[0,t-1]} | {\cal S}) \over P(q'_{[0,t-1]}|{\cal S})}  \label{sonD}\\
&& = P(q'_t| q'_{[0,t-1]}) \label{finalEqn}
\end{eqnarray}
where (\ref{defn00}) follows from Definition \ref{ClassAChannelDefinition}, (\ref{capD}) from the structure of a coding policy, and (\ref{sonD}) from the following inductive argument. Note that $P(q_0,q'_0|{\cal S}) = P(q_0,q'_0)$. If $P(q_{[0,t-1]},q'_{[0,t-1]} | {\cal S}) = P(q_{[0,t-1]},q'_{[0,t-1]})$, it follows that
\begin{eqnarray}
&& P(q_{[0,t]},q'_{[0,t]} | {\cal S})  \nonumber \\
&& = P(q'_t| q_{[0,t]}, q'_{[0,t-1]}, {\cal S}) P^*(q_t| q'_{[0,t-1]},q_{[0,t-1]}, {\cal S}) P(q_{[0,t-1]},q'_{[0,t-1]} | {\cal S})  \nonumber \\
&&= P(q'_t| q_{[0,t]}, q'_{[0,t-1]}) P^*(q_t| q'_{[0,t-1]},q_{[0,t-1]}) P(q_{[0,t-1]},q'_{[0,t-1]} | {\cal S})  \nonumber \\
&&= P(q'_t| q_{[0,t]}, q'_{[0,t-1]}) P^*(q_t| q'_{[0,t-1]},q_{[0,t-1]}) P(q_{[0,t-1]},q'_{[0,t-1]})  \nonumber \\
&& = P(q_{[0,t]},q'_{[0,t]})  
\end{eqnarray}
As a result, (\ref{sonD}) simplifies to (\ref{finalEqn}) by eliminating the conditioning on ${\cal S}$ and (\ref{EqCon}) holds.

We now use a similar argument as in (\ref{boundeNoAMS1}), but need to modify the steps due to the conditioning on ${\cal S}$:
%-----
\begin{eqnarray}
&&\lim_{T \to \infty} R_T \geq \limsup_{T \to \infty} {1 \over T} \bigg( \sum_{t=1}^{T-1} I(x_t; q'_t| q'_{[0,t-1]}, {\cal S}) + I(x_0;q'_0| {\cal S}) \bigg)\nonumber \\
 && = \limsup_{T \to \infty} {1 \over T} \sum_{t=1}^{T-1} I(x_t; q'_t| q'_{[0,t-1]}, {\cal S}) \nonumber \\
 &&=  \limsup_{T \to \infty} {1 \over T}  \sum_{t=1}^{T-1}  h_{\cal S}(x_t | q'_{[0,t-1]}) - h_{{\cal S}}(x_t|q'_{[0,t]}) \nonumber \\
 &&=  \limsup_{T \to \infty} {1 \over T}  \sum_{t=1}^{T-1}  h_{{\cal S}}(f(x_{t-1}) + B u_{t-1} +w_{t-1} | q'_{[0,t-1]}) - h_{{\cal S}}(x_t|q'_{[0,t]}) \nonumber \\
 && = \limsup_{T \to \infty} {1 \over T}  \sum_{t=1}^{T-1}  h_{{\cal S}}(f(x_{t-1}) + B u_{t-1} | q'_{[0,t-1]}) - h_{{\cal S}}(x_t|q'_{[0,t]}) \nonumber \\
 && = \limsup_{T \to \infty} {1 \over T}  \sum_{t=1}^{T-1} h_{{\cal S}}(f(x_{t-1})| q'_{[0,t-1]}) - h_{{\cal S}}(x_t|q'_{[0,t]}) \nonumber \\
 && = \limsup_{T \to \infty} {1 \over T}  \sum_{t=1}^{T-1} \bigg(E\bigg[P_{{\cal S}}(dx_{t-1}|q'_{[0,t-1]})  \log_2(|J(f(x_{t-1}))|) \bigg] \nonumber \\
 && \quad \quad \quad \quad \quad \quad \quad  \quad \quad \quad  + h_{{\cal S}}(x_{t-1} | q'_{[0,t-1]}) - h_{{\cal S}}(x_t|q'_{[0,t]})\bigg)  \label{indepW} \\
 &&\geq  V_{{\cal S}} - \liminf_{T \to \infty} \bigg( {1 \over T} h_{{\cal S}}(x_{T-1} | q'_{[0,T-1]}) \bigg)  \label{boundeNoAMS} 
\end{eqnarray}
 Here, 
 \begin{eqnarray}
 V_{{\cal S}} = \liminf_{T \to \infty} E\bigg[ {1 \over T}\bigg( \sum_{t=1}^{T-1} P_{{\cal S}}(dx_{t-1}|q'_{[0,t-1]}) \log_2(|J(f(x_{t-1}))|) \bigg) \bigg],\label{Vdefn}
 \end{eqnarray}
and (\ref{indepW}) follows from the fact that 
\begin{eqnarray} 
&& h_{{\cal S}}(f(x_{t-1}) | q'_{[0,t-1]}) = E\bigg[ \int P_{{\cal S}}(dx_{t-1}|q'_{[0,t-1]}) \log_2(|J(f(x_{t-1}))|)\bigg] + h_{{\cal S}}(x_{t-1} | q'_{[0,t-1]}),  \nonumber
\end{eqnarray} 
where the expectation is over the realizations of $q'_{[0,t-1]}$. Finally, we use the boundedness of $h(x_0)$ (and thus $h(x_0|q'_0)$) in (\ref{boundeNoAMS}). Thus, with $V_{{\cal S}} \geq L$, it follows that
\begin{eqnarray} \label{lowerBound1}
\liminf_{T \to \infty} \bigg( {1 \over T} h_{{\cal S}}(x_{T-1} | q'_{[0,T-1]}) \bigg) \geq L- C
\end{eqnarray}
We now seek to obtain an upper bound on $ h_{{\cal S}}(x_{T-1} | q'_{[0,T-1]})$. As in \cite{Matveev}, note that
\[ h_{{\cal S}}(x_{T}| q'_{[0,T]}) \leq h_{{\cal S}}(x_{T}, {\cal Y} | q'_{[0,T]}),\]
where ${\cal Y}$ is a binary random variable which is $1$ if $|x_{T}| \leq b(T)$ and $0$ otherwise. Let \[P_{\cal S}({\cal Y}=1) = P_{\cal S}(|x_{T}| \leq b(T)) =: p^{\cal S}_{T}.\]

Then,
\begin{eqnarray}\label{lowerBound2}
 h_{{\cal S}}(x_{T}, {\cal Y} | q'_{[0,T]}) &=& h_{{\cal S}}(x_{T} | q'_{[0,T]}, {\cal Y}) + H_{{\cal S}}({\cal Y} | q'_{[0,T]}) \nonumber \\
&\leq& h_{{\cal S}}(x_{T} | q'_{[0,T]}, {\cal Y}) + 1, \nonumber
\end{eqnarray}
since ${\cal Y}$ is binary. We have that
\begin{eqnarray}\label{lowerBound3}
&&h_{{\cal S}}(x_{T}| q'_{[0,T]}, {\cal Y}) \leq p^{\cal S}_{T}  {n \over 2} \log_2(2 \pi e b^2(T)) \nonumber \\
&& \quad \quad \quad \quad \quad \quad + (1 - p^{\cal S}_{T} ) h_{{\cal S}}\bigg(x_{T} \bigg| q'_{[0,T]}, |x_{T}| \geq b(T)\bigg) \nonumber
\end{eqnarray}
and
\begin{eqnarray}
&&h_{{\cal S}}\bigg(x_{T} \bigg| q'_{[0,T]}, |x_{T}| > b(T)\bigg) \nonumber \\
&&= h_{{\cal S}}\bigg(f(x_{T-1})+Bu_{T-1}+ w_{T-1}) \bigg| q'_{[0,T]}, |x_{T}| > b(T)\bigg) \nonumber \\
&& \leq h_{{\cal S}}\bigg(f(x_{T-1})+Bu_{T-1}+ w_{T-1}) \bigg| q'_{[0,T-1]},  |x_{T}| > b(T)\bigg) \label{uncondition1} \\
&&= h_{{\cal S}}\bigg(f(x_{T-1})+ w_{T-1}) \bigg| q'_{[0,T-1]}, |x_{T}| > b(T)\bigg) \nonumber \\
&&= h_{{\cal S}}\bigg(f(x_{T-1}) \bigg| q'_{[0,T-1]}, |x_{T}| > b(T)\bigg) \label{noiseOut} \\
&& =E\bigg[ \int P_{{\cal S} }\bigg(dx_{T-1} \bigg| q'_{[0,T-1]}, |x_{T}| > b(T)\bigg) \log_2\bigg(\bigg|J(f(x_{T-1}))\bigg|\bigg) \bigg] \nonumber \\
&& \quad \quad \quad \quad +  h_{{\cal S} }\bigg(x_{T-1} \bigg| q'_{[0,T-1]}, |x_{T}| > b(T)\bigg)\label{invertible} \\
%&& \leq E\bigg[ \int P_{{\cal S} }\bigg(dx_{T-1} \bigg| q'_{[0,T-1]}, |x_{T}| > b(T)\bigg) \log_2\bigg(\bigg|J(f(x_{T-1}))\bigg|\bigg) \bigg] \nonumber \\
%&& \quad \quad \quad \quad +  h_{{\cal S} }\bigg(x_{T-1} \bigg| q'_{[0,T-2]}, |x_{T}| > b(T)\bigg) \nonumber \\
&& \leq M +  h_{{\cal S} }\bigg(x_{T-1} \bigg| q'_{[0,T-1]}, |x_{T}| > b(T)\bigg) \nonumber \\
&& \quad \quad \vdots \nonumber \\
%&& \quad = h_{{\cal S}}\bigg(A^{T} (x_0 + \sum_{k=0}^{T-1} A^{-k-1} (w_k+u_k)) \bigg| q'_{[0,T]}, |x_{T}| > b(T) \bigg) \nonumber \\
&&  \leq M T + h_{\cal S}\bigg(x_0 \bigg|  |x_{T}| > b(T)\bigg)\label{lowerBound4} 
%&& \leq M T + {n \over 2} \log_2(2 \pi e K^2).\label{faydaliDenklem}
\end{eqnarray}
Here (\ref{uncondition1}) follows from that conditioning on a random variable reduces the differential entropy, and (\ref{noiseOut}) follows due to the fact that ${\cal S}$ determines the noise realizations. We note that non-linearity of $f$ add further technical issues when compared with the linear setup\footnote{Two technical intricacies here are as follows: For differential entropy (unlike discrete entropy) the relationship $h(x+y) \leq h(x) + h(y)$ does not in general hold for random variables $x, y$; this is why first a conditioning on ${\cal S}$ is taken in the proof. Furthermore, we cannot obtain an upper bound by taking out the conditioning on the event $|x_{T}| > b(T)$, since conditioning on a single event may decrease or increase entropy; note that conditioning on a random variable, however, does not increase the entropy.}. Here, $M$ is the supremum of $\log_2(|J(f(x))|)$.  In the above derivation in (\ref{invertible}), we use the fact that $f$ is invertible. In the last inequality, we use the fact that the entropy of a random variable with a fixed covariance is upper bounded by the entropy of a Gaussian with the same covariance, and that $|x_0|$ conditioned on ${\cal S}$ is upper bounded by $K^2$. 

Thus, by (\ref{lowerBound1}-\ref{lowerBound3}) and (\ref{lowerBound4})
we have
\begin{eqnarray}
&& \liminf_{T \to \infty} {1 \over T} \bigg( 1+ (1 - p^{\cal S}_{T}) \bigg(M T + h_{\cal S}(x_0 \bigg| |x_{T}| > b(T)) \bigg) +  p^{\cal S}_{T} {n \over 2} \log_2(2 \pi e b^2(T)) \bigg) \nonumber \\
&& \quad \quad  \quad \quad \quad \quad \quad \quad \geq L - C, \label{univBound}
\end{eqnarray}
Since $ h_{\cal S}\bigg(x_0 \bigg| |x_{T}| > b(T)\bigg) \leq (n/2) \log_2(2 \pi e K^2),$ it follows that for all $K$ and $\eta$
\[\limsup_{T \to \infty}P_{{\cal S}^K_{\eta}}(|x_{T}| \leq b(T)) \leq {M - (L-C) \over M},\]
for all $b(T)$ such that $\lim_{T \to \infty} \log_2(b(T))/T=0$. 
But now
\begin{eqnarray}
&& \limsup_{T \to \infty} P(|x_{T}| \leq b(T)) \nonumber \\
&& \leq \limsup_{T \to \infty} P(|x_{T}| \leq b(T), |x_0| \leq K) + \limsup_{T \to \infty} P(|x_{T}| \leq b(T), |x_0| \geq K) \nonumber \\
%&& = \limsup_{T \to \infty} P(|x_{T}| \leq b(T), |x_0| \leq K) + \limsup_{T \to \infty} P(|x_{T}| \leq b(T) | |x_0| \geq K) P(|x_0| \geq K) \nonumber \\
&& \leq \limsup_{T \to \infty} P(|x_{T}| \leq b(T), |x_0| \leq K) + P(|x_0| \geq K) \nonumber \\
&& = \limsup_{T \to \infty} \int P(d\eta)P_{{\cal S}^K_{\eta}}(|x_{T}| \leq b(T))  + P(|x_0| \geq K) \nonumber \\
&& \leq  \int P(d\eta) \limsup_{T \to \infty} P_{{\cal S}^K_{\eta}}(|x_{T}| \leq b(T))  + P(|x_0| \geq K) \label{FatouUsed} \\
&& \leq  \int P(d\eta) {M - (L-C) \over M} + P(|x_0| \geq K)\nonumber \\
&& = {M - (L-C) \over M} + P(|x_0| \geq K)\nonumber 
\end{eqnarray}
where we use Fatou's lemma in (\ref{FatouUsed}) and the fact that (\ref{univBound}) holds for every restriction of the noise realizations $\eta$ and $K$ values. Since an individual probability measure is tight, $\lim_{K \to \infty} P(|x_0| \geq K) = 0$, the right hand side can be made arbitrarily close to ${M - (L-C) \over M}$ and the result follows. \qed

\end{document}